\documentclass[12pt,a4paper]{article}
\usepackage[]{fontenc}
\usepackage[latin1]{inputenc}
 \usepackage[colorlinks]{hyperref}
\usepackage{a4,amsmath,amssymb,amsfonts,amsthm}
\usepackage{amsrefs}
\usepackage[all]{xy}

\def\mZ{\mathbb{Z}}
\def\mQ{\mathbb{Q}}
\def\mC{\mathbb{C}}

\def\mP{\mathbb{P}}

\def\Ev{{\rm Ev}}
\def\Spec{{\rm Spec}}

 \def\Tr{{\rm Tr}}

 \def\End{{\rm End}}
 \def\det{{{\rm det}}}
 
 \def\ch{{\rm ch}}
 
  \def\Spec{{\rm Spec}}
\def\Id{{\rm Id}}

\def\rk{{\rm rk}}

\def\depth{{\rm depth}}
\def\grade{{\rm grade}}
\def\Li{{\rm Li}}

\def\beginProof{\par{\bf Proof. }}
 \def\endProof{${\qed}$\par\smallskip}

 \def\Fp{{\mathfrak p}}
 \def\Fm{{\mathfrak m}}
\def\Fq{{\mathfrak q}}

 \def\mQ{{\Bbb Q}}
 
 \def\mZ{{\Bbb Z}}
 \def\mC{{\Bbb C}}
 
 \def\mF{{\Bbb F}}

 \def\CH{{\mathcal H}}
 \def\CO{{\mathcal O}}

 \def\CL{{\mathcal L}}

 \def\wt#1{\widetilde{#1}}
 
 \def\refeq#1{(\ref{#1})}

\def\P1{{{\bf P}^1}}
\def\cod{{\rm cod}}

\def\indlim{\underrightarrow{\lim}}

\def\Proj{{\rm Proj}}

\def\R{{\rm R}}

\def\L{{\rm L}}

\def\Sym{{\rm Sym}}
\def\indlim{\underrightarrow{\lim}}

\def\ch{{\rm ch}}

\def\max{{\rm max}}

\def\GTI{{\rm GTI}}
\def\cm{{\rm cm}}

\def\rTor{{\mathcal Tor}}

\def\reg{{\rm reg}}
\def\Rees{{\rm Rees}}
\def\height{{\rm height}}
\def\rad{{\rm radical}}
\def\TT{{\rm TT}}

\def\c1{{\rm c}_1}
\def\triv{{\rm triv}}
\def\r{{\rm r}}

 \newtheorem{theor}{Theorem}[section]
 \newtheorem{prop}[theor]{Proposition}

 \newtheorem{cor}[theor]{Corollary}
 \newtheorem{lemma}[theor]{Lemma}
 \newtheorem{sublem}[theor]{sub-lemma}

\newtheorem{rem}[theor]{Remark}
\newtheorem{ex}[theor]{Example}

\begin{document}

 \author{Damian R\"OSSLER\footnote{Mathematical Institute, 
University of Oxford, 
Andrew Wiles Building, 
Radcliffe Observatory Quarter, 
Woodstock Road, 
Oxford OX2 6GG, 
United Kingdom}}
 \title{The Riemann-Roch theorem in a singular setting}
\maketitle
\begin{abstract}
We prove a generalisation of the Grothendieck-Riemann-Roch theorem, which is valid for 
any proper and flat morphism between noetherian and separated schemes of odd characteristic. 
\end{abstract}

\tableofcontents

\section{Introduction}

The aim of this text is to prove a generalisation of the Grothendieck-Riemann-Roch (GRR) formula, which is valid for any proper and flat 
morphism of noetherian schemes. There is some speculation about such a generalisation in 
\cite[\S 2 \& \S 6.4]{SGA6}. 

We first recall a version of the GRR formula. 

 Write $K^0(W)$ for the Grothendieck group of vector bundles ($=$ coherent locally free sheaf) on a scheme $W$. The tensor product 
$\otimes$ endows this group with a canonical structure of commutative rings. There is also a ring endomorphism $\psi^2:K^0(W)\to K^0(W)$ (called the 2nd Adams operation), which sends 
a vector bundle $H$ to the element $\Sym^2(H)-\Lambda^2(H)$ (see \cite[I, \S 6]{FL-RR}). 

Let $S$ be a quasi-compact scheme, which carries an ample line bundle. 
Let $f:Y\to S$ be a perfect and proper morphism. The morphism $f$ then sends perfect 
complexes to strictly perfect complexes  since $S$ carries an  ample line bundle. 
Hence one may define a morphism of groups $\R^\bullet f_*:K^0(Y)\to K^0(S)$, which sends a vector bundle $V$ to the class of the strictly perfect complex 
$\R^\bullet f_*(V)$ in $K^0(S)$. 

For any vector bundle $H$, write 
$\Theta^2(H):=\bigoplus_i\Lambda^i(H)$. If $R$ is a commutative ring, 
$T$ is a $R$-algebra and $M$ is a $R$-module, write 
$M_T:=M\otimes_R T$. 

The GRR formula is then equivalent to the conjunction of the following statements.

(1) For any vector bundle $H$ on a scheme $W$, which carries an ample line bundle, the element $\Theta^2(H)$ is invertible in $K^0(W)_\mQ.$ 

(2) If $f$ is a Koszul-regular closed immersion, then for any vector bundle on $Y$ we have 
$$
\psi^2(\R^\bullet f_*(V))=\R^\bullet f_*(\Theta^2(N_f)\otimes\psi^2(V))
$$
in $K^0(S)_\mQ.$ Here $N_f$ is the conormal bundle of $f$. 

(3) If $f$ is smooth and projective then for any vector bundle on $Y$ we have 
\begin{equation}
\psi^2(\R^\bullet f_*(V))=\R^\bullet f_*(\Theta^2(\Omega_f)^{-1}\otimes\psi^2(V))
\label{eqARR}
\end{equation}
in $K^0(S)_\mQ.$ 

We note for future reference that in \refeq{eqARR}, we may have assumed without restriction of generality that the sheaves $\R^i f_*(V)$ are locally free, because 
any vector bundle on $Y$ admits a resolution by vector bundles with this property.

We refer to \cite[VIII]{SGA6} and \cite[V, \S 7]{FL-RR} for details (see also Proposition \ref{propLCIT} below). The equalities (2) and (3) can be joined to give a Riemann-Roch formula for any $f$ which admits a factorisation into 
a Koszul-regular closed immersion followed by a smooth and projective morphism (or in other words, for any $f$ which is projective and lci). This version of the GRR formula 
is often called the Adams-Riemann-Roch (ARR) theorem for the Adams operation $\psi^2$.

Our aim is to formulate a generalisation of (3), which will work without assumption of smoothness or projectivity for 
$f$ but only retains the assumption of flatness and of properness.  Note that one expects a priori that 
any generalisation of (3) to the proper and flat case must be formulated using the Grothendieck group of coherent sheaves rather than the Grothendieck group of locally free coherent sheaves, because 
the sheaves of differentials of non smooth morphism are not be perfect complexes in general. 
One can also seek to generalise (2) but this would presumably require techniques different from the ones 
that will be employed in this article. 

We shall now formulate the generalisation of (3), which is the main result of this article. We first need to introduce a few notions and some terminology. 
 
If $W$ is a scheme, we shall write $K_0(W)$ for the Grothendieck group of coherent sheaves on $W$. 
The group $K_0(W)$ is a $K^0(W)$-module via the tensor product of coherent sheaves by locally free sheaves. We let $\r K^0(W)$ be the quotient of 
$K^0(W)$ by the annihilator of $K_0(Y)$. The group $K_0(W)$ then obtains a $\r K^0(W)$-module structure. Note also that the natural map $K^0(W) \to K_0(W)$ factors 
through $\r K^0(W)$, so we can speak of the image of an element of $\r K^0(W)$ in $K_0(W)$. 
Similarly, we will write $\r K^0(W)_\mQ$ for the quotient of $K^0(W)_\mQ$ by the annihilator 
of $K_0(W)_\mQ$. Note that there is a natural map of $\mQ$-vector spaces 
$(\r K^0(W))_\mQ\to \r K^0(W)_\mQ$ but it is not clear that this map is an isomorphism in general. 

Suppose from now on that $f$ is flat and proper (no other assumptions). We shall also suppose that $S$ is a separated noetherian scheme and that $2$ is invertible on $S$ (it will become apparent 
below why this is necessary). We do not assume that $S$ carries an ample line bundle anymore.

Let $\Delta\subseteq Y\times_S Y$ be the relative diagonal of $Y\to S$. Let $I_\Delta\subseteq\CO_{Y\times_S Y}$ 
 be the sheaf of ideals of $\Delta.$

Let $\pi:\wt{X}\to Y\times_S Y$ be 
the blow-up of $Y\times_S Y$ along $\Delta$ and let $\phi:E\to\Delta$ be the corresponding 
exceptional divisor. Let $N_{E/\wt{X}}$
 be the conormal bundle of $E$ in $\wt{X}.$ Let \mbox{$\cm(f)=\cm(Y\to S)\geq 0$} be the minimal natural number $\lambda$ such that 
$\R^a\pi_*(\CO(-E)^{\otimes r})=0$ for all $a>0$ and such that the natural 
morphism of sheaves $I_\Delta^r\to \pi_*(\CO(-E)^{\otimes r})$ is an isomorphism  for all $r\geq \lambda.$ This exists by \cite[Cohomology of Schemes, Lemmas 14.2 and 14.3]{StacksProject}. Note that if 
$Y$ is smooth over $S$ then  $\cm(f)=0$. Note also that the invariant $\cm(Y\to S)$ makes 
sense for any scheme $Y$, which is separated and of finite type over $S$. The integer $\cm(f)$ is bounded by the Castelnuovo-Mumford regularity of 
the Rees Algebra $\oplus_{k\geq 0}I_\Delta^k$. See section \ref{secCM} below. 

Write  
$s(l,j)$ for the Stirling numbers of the first kind and $E_j$ for the $j$-th Euler number. 
By definition
$$
t(t-1)\dots(t-l+1)=\sum_{j} s(l,j)t^j
$$
(where it is understood that $t(t-1)\cdots(t-k+1)=1$ if $k=0$) 
and 
$$
{2\over e^t+1}=\sum_{j=0}^\infty E_j {t^j\over j!}.
$$
Note that for $j\geq 1$ we have 
$$E_j=2(-1)^{j+1}(2^{1+j}-1)\zeta_\mQ(-j)=2(-1)^{j}(2^{1+j}-1){B_{j+1}\over (j+1)}$$
where the $B_{j+1}$ are Bernoulli numbers and $\zeta_\mQ(\cdot)$ is the Riemann z\^eta function (see 
\cite[chap. IV]{Washington-Cyc}).

We will show below that for any line bundle $\CL$ on a noetherian scheme $W$, the  element $\CL-1$  is nilpotent in $\r K^0(W)_\mQ$ (in fact even in $\r K^0(W)$) and we shall write $\delta_0(\CL)$ for the smallest natural number $n$ such that $(\CL-1)^{\otimes (n+1)}=0$ in $\r K^0(W)_\mQ$. This simple fact is what will allow us to circumvent the projectivity hypothesis on $f$. Now let
$$\delta_0(f):=\delta_0(N_{E/\wt{X}})$$

and for any $\delta,\lambda\geq 0$ let
$$
\GTI(f,\delta,\lambda):={(-1)^\lambda\over 2}\sum_{j=0}^\delta E_j \sum_{k=0}^{\delta}\sum_{u=0}^k\ {(-1)^{k-u}s(k,j)\over u!(k-u)!}I^{u+\lambda}_\Delta/I^{u+\lambda+1}_\Delta+\sum_{k=0}^{\lambda-1}(-1)^k I^k_\Delta/I^{k+1}_\Delta\in K_0(Y)_\mQ.
$$ 
We will show in section \ref{secMT} below that $\GTI(f,\delta,\lambda)$ is constant in the range $\delta\geq\delta_0(f)$ and $\lambda\geq\cm(f)$. 
We shall write $\GTI(f)$ for this constant.

We will show that if $Y$ has an ample line bundle and is of finite 
dimension, then \mbox{$\delta_0(f)\leq\dim(E)$} (see Lemma \ref{lemPOL} (c) (3)). 
In particular, we then have 
$$\delta_0(f)\leq\max\{\dim(Y_s)\,|\,s\in S\}+\dim(Y)$$ 
(use eg \cite[4.3.12]{Liu-Alg} and \cite[Th. 15.17]{Mat-Co}) and if $S$ is the spectrum of a field and 
$Y$ is integral and projective over $S$, we even have $\delta_0(f)\leq 2\dim(Y)-1$. 

\begin{theor} Suppose that $V$ is a vector bundle on $Y$ and 
that $\R^i f_*(V)$ is locally free for all $i\geq 0$. Then 
the equality 
\begin{eqnarray}
&&\sum_i(-1)^i[\Sym^2(\R^i f_*(V)-\Lambda^2(\R^i f_*(V))]=\R^\bullet f_*(\GTI(f)\otimes(\Sym^2(V)-\Lambda^2(V)))
\label{eqFForm}
\end{eqnarray}
holds in $K_0(S)_\mQ.$
\label{theorMT}
\end{theor}

We will show in Corollary \ref{lemMIR} below that for any vector bundle $V$ on a noetherian scheme $W$ with an ample family of line bundles, we have 
\begin{equation}
\Theta^2(V)^{-1}={1\over 2}\sum_{j=0}^\delta E_j \sum_{k=0}^{\delta}\sum_{u=0}^k\ {(-1)^{k-u}s(k,j)\over u!(k-u)!}\Sym^u(V)
\label{eqRDE}
\end{equation}
in $K^0(W)_\mQ$ for any $\delta$ such that $(\CO(1)-1)^{\otimes(\delta+1)}=0$ in $K^0(P)_\mQ$. Here 
 $P:=\Proj(\Sym^\bullet(V))$ and $\CO(1)$ is the tautological line bundle on $P.$ 
 As a consequence of \refeq{eqRDE}, we see that \refeq{eqFForm} naturally reduces to \refeq{eqARR} when $f$ is smooth and projective.

\begin{cor} Suppose that $S$ is the spectrum of a field and that $Y$ is integral and projective over $S$. Let $V$ be a vector bundle on $Y$. 

{\rm (a)} If $\dim(Y)\leq 1$ we have
\begin{eqnarray*}
&&(-1)^\lambda\chi(Y,V)\\&=&
\chi(Y,\big[{3\over 4}I^{\lambda}_\Delta/I^{1+\lambda}_\Delta-
{1\over 4}I^{1+\lambda}_\Delta/I^{2+\lambda}_\Delta+(-1)^\lambda\sum_{k=0}^{\lambda-1}(-1)^k I^k_\Delta/I^{k+1}_\Delta\big]\otimes
\big(\Sym^2(V)-\Lambda^2(V)\big))
\end{eqnarray*}
for all $\lambda\geq\cm(f)$. 

{\rm (b)} If $\dim(Y)\leq 2$ we have
\begin{eqnarray*}
&&(-1)^\lambda\chi(Y,V)\\&=&
\chi(Y,\big[{15\over 16}I^{\lambda}_\Delta/I^{1+\lambda}_\Delta-
{11\over 16}I^{1+\lambda}_\Delta/I^{2+\lambda}_\Delta
+{5\over 16}I^{2+\lambda}_\Delta/I^{3+\lambda}_\Delta
-{1\over 16}I^{3+\lambda}_\Delta/I^{4+\lambda}_\Delta
+(-1)^\lambda\sum_{k=0}^{\lambda-1}(-1)^k I^k_\Delta/I^{k+1}_\Delta\big]\\&\otimes&
\big(\Sym^2(V)-\Lambda^2(V)\big))
\end{eqnarray*}
for all $\lambda\geq\cm(f)$. 
\label{corRRV}
\end{cor}
Here $\chi(Y,\cdot)$ takes the Euler characteristic of a coherent sheaf.

\begin{rem}\rm We note the following conceptual consequence of Theorem \ref{theorMT}. The correction factor $\GTI(f)$ 
depends only on the thickening of order $\cm(f)+\delta_0(f)$ of the relative diagonal of 
$Y\to S$ (this thickening is the algebra of differential operators of order $\cm(f)+\delta_0(f)$ of $Y\to S$;  see \cite[IV.4, 16.7]{EGA}). 
As we have seen, the invariant $\delta_0(f)$ is usually easy to estimate  but by contrast $\cm(f)$ depends on the singularities of the fibres of the morphism $f$ and is in general difficult to compute (although it is in principle effectively computable in any given case). One expects 
$\cm(f)$ to be large if the fibres of $Y$ are 'very' singular. In particular, one expects a Riemann-Roch 
theorem for a variety with 'complicated' singularities to have a correction factor which depends on 
a big infinitesimal neighbourhood of the diagonal.\end{rem}

\begin{rem}\rm It is instructive to carry out the comparison between formula (b) in Corollary \ref{corRRV}   and the output of the classical Riemann-Roch theorem for smooth projective surfaces (to convince the reader that formula (b) really
does generalise that theorem to the singular setting\dots). 

So suppose for the time of the present paragraph that $Y$ is a smooth projective surface. Let $D$ be a divisor on $Y$.  Then formula (b) applied to $V=\CO(D)$ 
and $V=\CO_Y$ separately implies that we have
\begin{eqnarray*}
&&\chi(Y,\CO(D)-\CO_Y)\\&=&\chi(Y,\big({15\over 16}\CO_Y-{11\over 16}\Omega_Y+{5\over 16}\Sym^2(\Omega_Y)-{1\over 16}\Sym^3(\Omega_Y)\big)\otimes(\CO(D)^{\otimes 2}-\CO_Y)).
\end{eqnarray*}
Now for any divisors $H$ and $J$ on $Y$, write 
$\langle H,J\rangle$ for the intersection number of $H$ and $J$. From the definition of this pairing (see \cite[V.1]{Hartshorne-AG}), we see that we have 
$$\chi(Y,(\CO(H)-\CO_Y)\otimes(\CO(J)-\CO_Y))=\langle H,J\rangle.$$ We will also use the fact that for any two vector bundles $V_1$ and $V_2$ on $Y$, we have 
$$
\chi(Y,(V_1-V_2)\otimes(\CO(H)-\CO(J))=0
$$ if $\rk(V_1)=\rk(V_2)$ and $\det(V_1)\simeq\det(V_2).$ This follows from the fact that the Chern character of $(V_1-V_2)\otimes(\CO(H)-\CO(J))$ in the (Chow) intersection ring of $Y$ vanishes in this situation. Now using the splitting principle and the theorem on symmetric functions, we may compute that
\begin{eqnarray*}
&&{15\over 16}\CO_Y-{11\over 16}\Omega_Y+{5\over 16}\Sym^2(\Omega_Y)-{1\over 16}\Sym^3(\Omega_Y)\\&=&
{1\over 4} - {1\over 8}(\Omega_Y-2\CO_Y) - {1\over 16}(\det(\Omega_Y)-\Omega_Y+\CO_Y)\\&+& {1\over 8}((\det(\Omega_Y)-\Omega_Y+\CO_Y)\otimes(\Omega_Y-2\CO_Y)) + 
  {1\over 16}(\Omega_Y-2\CO_Y)^{\otimes 2} - {1\over 16}(\Omega_Y-2\CO_Y)^{\otimes 3}
  \end{eqnarray*}
in $K^0(Y)_\mQ$ (this follows from the polynomial identity 
\begin{eqnarray*}
&&{15\over 16} - {11\over 16}(x + y) + {5\over 16}(x^2 + xy + y^2) - {1\over 16}(x^3 + 
     x^2y + y^2x + y^3)\\&=&
     {1\over 4} - {1\over 8}(x+y-2) - {1\over 16}(x-1)(y-1) \\&+& {1\over 8}((x-1)(y-1)(x+y-2)) + 
  {1\over 16}(x+y-2)^2 - {1\over 16}(x+y-2)^3\,\,\,\,).
\end{eqnarray*}
  
and thus, using the just mentioned computational rules, we have 
\begin{eqnarray*}
&&\chi(Y,\CO(D)-\CO_Y)\\&=&\chi(Y,\big({15\over 16}\CO_Y-{11\over 16}\Omega_Y+{5\over 16}\Sym^2(\Omega_Y)-{1\over 16}\Sym^3(\Omega_Y)\big)\otimes(\CO(D)^{\otimes 2}-\CO_Y))\\
&=&
\chi(Y,({1\over 4} - {1\over 8}(\Omega_Y-2\CO_Y))\otimes (\CO(D)^{\otimes 2}-\CO_Y)))=
{1\over 4}\chi(Y,\CO(D)^{\otimes 2}-\CO_Y)\\&-&{1\over 8}\chi(Y,\big[\Omega_Y-2\CO_Y+\det(\Omega_Y)-\Omega_Y+\CO_Y\big]\otimes (\CO(D)^{\otimes 2}-\CO_Y)))
\\&=&
{1\over 4}\chi(Y,\CO(D)^{\otimes 2}-\CO_Y)-{1\over 8}\chi(Y,(\det(\Omega_Y)-\CO_Y)\otimes (\CO(D)^{\otimes 2}-\CO_Y)))\\&=&
{1\over 4}\langle D,D\rangle+{2\over 4}\chi(Y,\CO(D))-{2\over 4}\chi(Y,\CO_Y)-{1\over 8}\langle K,2D\rangle\\
&=&
{1\over 2}\langle D,D\rangle -{1\over 2}\langle K,D\rangle
\end{eqnarray*}
where $K$ is a divisor representing the line bundle $\det(\Omega_Y).$ This is the classical Riemann-Roch formula for surfaces \cite[V.1, Th. 1.6]{Hartshorne-AG}.\end{rem} 

We shall now consider some examples. Remember the running assumptions on $Y\to S$: the morphism $Y\to S$ is flat and proper,  $S$ is separated and noetherian, and $2$ is invertible on $S$. The properties of $Y\to S$ considered in the examples are all in addition to these properties. 

\begin{ex}\rm When $f$ is smooth then $\cm(f)=0$ and we have 
$$
\GTI(f)={1\over 2}\sum_{j=0}^\delta E_j \sum_{k=0}^{\delta}\sum_{u=0}^k\ {(-1)^{k-u}s(k,j)\over u!(k-u)!}\Sym^u(V)
$$
for any $\delta\geq\delta_0(f).$ Note however that in this situation, the element 
$\Theta^2(\Omega_f)$ is not necessarily invertible in $K^0(Y)_\mQ$ and 
$\GTI(f)$ cannot a priori be compared with $\Theta^2(\Omega_f)$ (an ample family of line bundles is needed for this).\end{ex}

\begin{ex}\rm We will show in Proposition \ref{propSLCI} below that if $S$ is an integral Cohen-Macaulay scheme whose regular locus contains an open set and whose local rings have infinite residue fields and if $Y\to S$ has geometrically integral Cohen-Macaulay fibres with at most hypersurface singularities (see before \ref{propSLCI} for the definition), then $\cm(f)=0$. In addition, we will show that 
in the situation of the last sentence, we  have 
$$
I^k_\Delta/I^{k+1}_\Delta\simeq \Sym^k(\Omega_f)
$$
although $\Omega_f$ might not be locally free. 
These conditions are met in particular in the situation where $S$ is an integral Cohen-Macaulay scheme whose local rings have infinite residue fields and $Y$ can be embedded locally on $Y$ (as a $S$-scheme) as a Cartier divisor into a smooth scheme over $S$. If there is a (global) embedding of $Y$ as a Cartier divisor into a smooth scheme over $S$, then the GRR theorem for 
lci morphisms and again one can show that Theorem \ref{theorMT} reduces 
to GRR in that case, although the argument for showing this is quite indirect (see Proposition \ref{propLCIT} below). See 
\cite{Savin-MO} for an example of a scheme, which satisfies the just described conditions but does not afford a global embedding into a smooth scheme.\end{ex}

\begin{ex}\rm At the other end of the spectrum, one might consider finite, flat and purely inseparable 
morphisms. Here is a simple example. 
Suppose that $Y$ is a smooth curve over a noetherian separated scheme $T$. Suppose that $T$ is a scheme of characteristic $p>0$, where $p$ is an odd prime. 
Let $F_{Y/T}:Y\to Y^{(p)}$ be the relative Frobenius morphism and set $f:=F_{Y/S}$ and $S:=Y^{(p)}$. 
The ideal $I_\Delta$ is then nilpotent, $E$ is empty and $f$ is finite (and flat). Suppose for simplicity that $V=\CO_Y.$ Theorem \ref{theorMT} then states that 
\begin{eqnarray*}
&&\Sym^2(f_*(\CO_Y))-\Lambda^2(f_*(\CO_Y))=
\sum_{k=0}^{\lambda-1}(-1)^k f_*(I^k_\Delta/I^{k+1}_\Delta)
\end{eqnarray*}
in $K_0(S)_\mQ$ if $\lambda\geq\cm(f).$ From the definitions, we see that the integer $\cm(f)$ is the smallest 
natural number $\lambda$ such that $I_\Delta^\lambda=0$. 
It is shown in \cite[\S 2]{PR-ARR} that $\cm(f)=p$ and that $I^k_\Delta/I^{k+1}_\Delta\simeq\Omega_{Y/T}^{\otimes k}$. We thus have an equality
\begin{eqnarray*}
&&\Sym^2(f_*(\CO_Y))-\Lambda^2(f_*(\CO_Y))=
\sum_{k=0}^{p-1}(-1)^k f_*(\Omega^{\otimes k}_{Y/T})
\end{eqnarray*}
in $K_0(S)_\mQ$. 
On the other hand, the morphism $f$ is lci because it is a $T$-morphism between two smooth schemes over $T$. The Adams-Riemann-Roch for the operation $\psi^2$ thus applies and gives the identity
\begin{eqnarray*}
\hskip-0.5cm
&&\Sym^2(f_*(\CO_Y))-\Lambda^2(f_*(\CO_Y))=
f_*(f^*(1+\Omega_{S/T})\otimes(1+\Omega_{Y/T})^{-1})=
f_*((1+\Omega_{Y/T}^{\otimes p})\otimes(1+\Omega_{Y/T})^{-1})
\end{eqnarray*}
in $K^0(S)_\mQ$ (in fact the equality holds in $K^0(S)[{1\over 2}]$). Now, since 
$\Omega_{Y/T}$ is a line bundle, we have 
$$
(1+\Omega_{Y/T}^{\otimes p})\otimes(1+\Omega_{Y/T})^{-1}=\sum_{k=0}^{p-1}(-1)^k \Omega^{\otimes k}_{Y/T}
$$
so Theorem \ref{theorMT} gives the same formula as ARR in this situation (but Theorem \ref{theorMT} only provides an identity in the Grothendieck group of coherent sheaves). We will see that the proof of Theorem \ref{theorMT} is completely elementary in this situation, because it does not involve any blow-up construction and only relies on the nilpotence of 
$I_\Delta$ and its equivariance under the natural involution of $X\times_S X$. It is curious that this elementary 
argument was never discovered before. The argument above can be generalised to 
higher relative dimension (over $T$). We leave the details as an exercise for the interested reader.\end{ex} 

\begin{ex}\rm  Here is a numerical example in a situation where $\cm(f)\not=0$ and $f$ is not finite or lci. 
We consider the plane projective curve $Y$ over $S=\Spec(\mF_5)$ defined by the equations $zx,z^3$. This is a non-reduced curve 
whose underlying reduced scheme is a copy of the projective line. It carries a "thickening" of order $3$ of the origin and it is not lci (this can be shown directly but it also follows from the numerical 
calculations below together with Proposition \ref{propSLCI}). In this case, the blow-up morphism of the diagonal of $Y\times_S Y$ is finite over $Y\times_S Y$ and we thus have $\delta_0(f)\leq 1$. We would like to compute
$\chi(Y, \GTI(f,1,\lambda))$ for $\lambda\geq 0$. For this it is sufficient to be able to compute 
the quantity $\chi(Y\times_S Y,\CO_{Y\times_S Y}/I_\Delta^n)$ for $1\leq n\leq \lambda+2$ (by the formula in 
Corollary \ref{corRRV} (a), which is valid because $\delta_0(f)\leq 1$). We made use of the 
computer package Magma (see \cite{Magma}) to do this for small values of $n.$ This computer package 
has a routine, which computes the dimension of the cohomology groups of the structure sheaf of 
a projective scheme over a field described by homogenous equations; the underlying algorithm is based on the Beilinson-Gelfand-Gelfand combinatorial description of the derived category of projective space (see \cite{BGG}). 
We first embed $\mP^2_S\times_S\mP^2_S$ into $\mP^8_S$ via the Segre embedding and we 
find explicit homogenous equations in $\mP^8_S$ for the diagonal of $\mP^2_S\times_S\mP^2_S$ and for $Y\times_S Y\subseteq \mP^2_S\times_S\mP^2_S$ . We can then work entirely in $\mP^8_S$ and we can use the command {\it CohomologyDimension} of 
Magma to  compute that 
$$
\chi(Y\times_S Y,\CO_{Y\times_S Y}/I_\Delta^n)=3, 5, 7, 7, 5, 0, -7, -16, -27, -40
$$
for $n$ running from $1$ to $10$ (the computation took about two hours on the computer mainframe of 
the Oxford Mathematical Institute). 
From this, we can compute 
$$
\chi(Y,\GTI(f,1,\lambda))={7\over 4},2, {5\over 2}, {5\over 2}, {11\over 4}, 3, 3, 3, 3
$$
for $\lambda$ running from $0$ to $8$. So $\chi(Y,\GTI(f,1,\lambda))$ has the 
value $\chi(Y,\CO_Y)=3$ required by Corollary \ref{corRRV} (a) for $5\leq\lambda\leq 8$. 
This suggests (but does not prove) that $\cm(f)\leq 8$ in the present situation. In particular, 
it suggests that $\chi(Y,\GTI(f,1,\lambda))=\chi(Y,\GTI(f)=3$ for all $\lambda\geq 5.$
\end{ex}

\begin{ex}\rm Theorem \ref{theorMT} implies that the Euler characteristic of certain linear combinations 
of the sheaves $I^r_\Delta/I^{r+1}_\Delta$ are divisible by certain integers. Eg when $Y$ is integral of 
dimension $\leq 2$ and 
$S$ is a field, Corollary \ref{corRRV} implies that 
\begin{equation}
15\chi(Y,I^{\lambda}_\Delta/I^{1+\lambda}_\Delta)-
{11}\chi(Y,I^{1+\lambda}_\Delta/I^{2+\lambda}_\Delta)
+{5}\chi(Y,I^{2+\lambda}_\Delta/I^{3+\lambda}_\Delta)
-\chi(Y,I^{3+\lambda}_\Delta/I^{4+\lambda}_\Delta)
\label{eqDIVP}
\end{equation}
is divisible by $16$ if $\lambda\geq\cm(f).$ In particular, this gives an obstruction for 
$\cm(f)$ to vanish in that situation. Similar divisibility properties can be worked out in any dimension. It would 
be interesting to have a general formula. Note that a direct proof of the fact that \refeq{eqDIVP} is divisible by $16$  when $Y$ is a smooth surface over $S$ (so that  $I^r_\Delta/I^{r+1}_\Delta\simeq\Sym^r(\Omega_Y)$) involves long winded computations with symmetric functions.\end{ex}

{\bf Outline of the proof.} Our method of proof is based on a fundamental remark of Nori, who noticed in \cite{Nori-HRR}
that the Adams-Riemann-Roch for a smooth variety over a field can be understood as a special case of the geometric fixed point formula of Atiyah-Bott. 
More precisely, he showed that the Adams-Riemann-Roch theorem (for $\psi^2$) for a smooth and projective variety $Y$  over a field of odd characteristic can be obtained 
by applying the fixed point formula of Atiyah-Bott to the involution of $Y\times_S Y$, which swaps the factors. 
We generalise this method to the singular and relative setting, using instead a more general fixed point formula 
proven by Thomason. We show that the local term of Thomason's fixed point formula for  a diagonalisable group of order $2$ can be explicitly computed by blowing-up the fixed point set and analysing the asymptotic behaviour of the local terms associated with increasing powers of the ideal of the exceptional divisor. As mentioned above, an essential point here is that the action of 
a line bundle minus the unit on the Grothendieck group of coherent sheaves of a noetherian scheme is nilpotent. This is easy to prove but it allows us to forego any projectivity assumptions. Once the 
local term of a sufficiently high power $n$ of the ideal of the exceptional divisor is computed, we can 
push back down to $Y\times_S Y$. The size of $n$ is controlled by the Castelnuovo-Mumford regularity 
of the Rees algebra of the ideal of the diagonal on $Y\times_S Y$ and this leads to the integer $\cm(f)$, which is in essence a measure of the singularity of the fibres of $Y\to S$. Note that our computations are greatly simplified by the fact that the group action is of order $2$. It is presumably 
also possible to provide explicit formulae for the local term of group actions of higher order but 
for that one would probably have to consider several successive blow-ups (about this, see 
\cite[p. 278]{Fresan-Periods}), leading to combinatorial problems that do not arise for group actions of order $2$. The restriction to odd characteristic in Theorem \ref{theorMT} comes from the fact that to apply Thomason's formula, one needs to know that the group $\mZ/2\mZ$ is diagonalisable, and this fails in characteristic $2.$ It seems unlikely that a singular generalisation of the GRR exists, which 
is sufficiently explicit and avoids any restriction on the characteristic. For explanations about this, see remark \ref{remRC} below.

{\bf Structure of the article.} The structure of the article is a follows. In section \ref{secEQSEC}, we prove equation \refeq{eqRDE}. This does not involve the fixed point formula and involves 
a limiting process and some complex analysis. In section \ref{secFP}, we review Thomason's fixed point 
formula and we give an explicit formula for the local term, by blowing-up the fixed point scheme. In section \ref{secREG}, we provide a different computation of the local term, in the situation where the scheme under scrutiny can be equivariantly embedded by a closed immersion of finite tor-dimension into a scheme, whose fixed point scheme is regularly embedded. This was in essence already done by Thomason in \cite[Th. 3.5]{Thom-L} but his result 
is not general enough for our purposes. The results of this section are used in section 
\ref{secSC} to show that the ARR formula and the formula of Theorem \ref{theorMT} coincide when the morphism is lci. In section \ref{secCM}, we review the definition of the Castelnuovo-Mumford regularity of a graded ring and we collect some results  on the Castelnuovo-Mumford regularity of Rees algebras, which are available in the literature. The computation of the Castelnuovo-Mumford regularity of a Rees algebra is a problem studied by several people since the late 1970s, in particular Vasconcelos and his school. 
In section \ref{secMT}, we prove Theorem \ref{theorMT} by applying the results 
of section \ref{secFP} to the involution of $Y\times_S Y$ swapping the factors. 
In section \ref{secSC}, we consider Theorem \ref{theorMT} in the situation where 
$Y\to S$ is lci, and also in the situation where the fibres of $Y\to S$ have Cohen-Macaulay  hypersurface singularities and are geometrically integral. 
In the latter situation, Theorem \ref{theorMT} formally looks very much like the ARR theorem in the smooth case. The 
crucial point here is that the diagonal immersion of a geometrically integral scheme of finite type over a field with Cohen-Macaulay  hypersurface singularities 
is an almost complete intersection in the sense of \cite{HSV-ACI}, and the sheaf of ideals of the diagonal is then locally generated by a $d$-sequence in the sense of Huneke (see \cite{Huneke-d}); this allows us to 
show that $\cm(f)=0$ using the results collected in section \ref{secCM}.

\begin{rem} \rm Note that there is already a singular form of the GRR 
theorem in the literature, namely the singular Riemann-Roch theorem of Baum-Fulton-MacPherson (see \cite{BFM-RR} - note that this article only treats varieties over fields, but the method could presumably be generalised). 
This theorem is formulated in terms of Chow homology but it could presumably 
be translated into a formula involving only Euler characteristics. The conceptual difference between the theorem of Baum-Fulton-MacPherson and Theorem \ref{theorMT} is that 
the right-hand side of the formula proven by these authors is computed using immersions into 
smooth schemes, whereas in Theorem \ref{theorMT} the right-hand side is computed directly on $Y$. There is a similar state of affairs in the theory of Grothendieck duality. The dualising 
complex of a projective scheme can be computed using an ambient projective 
space (as in \cite[III.7]{Hartshorne-AG}) or it can be described intrinsically using residual complexes (see \cite{Hartshorne-Residues}, or 
\cite[II.7, Th. 7.14.2]{Hartshorne-Residues} for the case of curves). In this sense, Theorem \ref{theorMT} provides the analogue for the Riemann-Roch theorem of the intrinsic description of the dualising complex in Grothendieck duality. Note that we also provide in section \ref{secREG} a computation of the right-hand side of the Riemann-Roch formula in the situation where an embedding into a smooth scheme is available.\end{rem}

\begin{rem} \rm The attentive reader might have noticed that we have outlined two different 
proofs 
of the equivalence between the ARR formula and the formula of Theorem \ref{theorMT}. 
The first one, given in section \ref{secEQSEC}, does not involve any reference 
to the fixed point formula, and works only when $Y\to S$ is smooth. 
The second one, given in section \ref{secSC}, works for all lci morphisms. However, this second proof relies on the 
unicity of the local term of the fixed point formula, which implies that the local terms computed in section \ref{secREG} and in section \ref{secFP} coincide when $Y\to S$ is lci. 
It is desirable to find a way to compare the two formulae directly when $Y\to S$ is lci, 
without resorting to this unicity. In other words, one would like to have a direct combinatorial proof of the equivalence of the formulae. 
The results of  \cite[Th. 6.3 and Prop. 10.3]{Quillen-Homology} should be relevant here  but we don't how to apply them.
A related question  is: is there a simple upper bound for 
$\cm(Y\to S)$ when $Y\to S$ is lci? Also, is it true that, when $Y$ is lci over $S$,  
the ideal of the diagonal can be locally generated by a sequence of regular type $r$ (in the 
sense of \cite[before Th. 1.3]{Trung-C}) for some $r$? Note that in section \ref{secSC} we answer all these questions in the 
situation where $Y\to S$ factors through a smooth scheme as a Cartier divisor (but our answer 
does not rely on  \cite{Quillen-Homology}).\end{rem}

\begin{rem}\rm Suppose $S$ is a field. It would be interesting to compute $\cm(Y\to S)$ in terms of the structure of 
the singular points of $Y$ when $Y$ has isolated singularities. Can any bounds for $\cm(Y\to S)$  
be given for certain classes of singularities (eg rational singularities)? \end{rem}

\begin{rem}\rm Each Adams operation $\psi^k$ ($k\geq 2$) has an associated Adams-Riemann-Roch theorem (see \cite[V.7]{FL-RR}).  These Riemann-Roch theorems  are all equivalent to 
the GRR formula, at least if one works with coefficients in $\mQ$. So at first sight it does not seem very interesting to model  a singular generalisation of the GRR formula on Adams operations $\psi^k$ for 
natural numbers $k>2$. 
However, by doing so, one could presumably avoid the restriction to odd characteristic (excluding the characteristics prime to $k$ instead). To carry this out, one would probably have to consider the fixed point formula for cyclic permutations of the  
$k$-fold product $Y\times_S Y\times_S\cdots\times_S Y$ (this is done in \cite{Nori-HRR} in the smooth case and when $S$ is the spectrum of field). As explained above, 
this would lead to combinatorial problems that we can forego in the case $\psi^k=\psi^2.$\end{rem}

\begin{rem}\rm\label{remRC} If one wanted to drop the restriction on the characteristic, one would 
have to consider the derived functors (in the sense of Quillen) of the non additive functor 
$\Sym^2(\cdot)$. One can forego having to consider these derived functors when one deals only with 
locally free sheaves (because the higher derived functors of $\Sym^2$ of a free module vanish). This 
is why no restriction on the characteristic is made in the classical ARR theorem, where all the sheaves in sight are locally free, including in the cotangent complex, which can be locally represented by a complex of locally free sheaves when the morphism is lci. When one considers non lci morphisms, coherent sheaves appear and without any restriction on the characteristic, consideration of the derived functors of 
$\Sym^2(\cdot)$ is inevitable. It would be interesting to see such calculations (which would in particular 
involve revisiting Thomason's fixed point formula) but the resulting formula would likely be very complicated. About the ARR theorem and the derived functors of 
$\Sym^2(\cdot)$, see also \cite{Koeck-KC}.\end{rem}

\section{Comparison between old and new}

\label{secEQSEC}
The aim of this section is mainly to prove equality \refeq{eqRDE}.

We need a few preliminary results. We will use the terminology described in the introduction.

We define
$$
\TT(x,t):=\sum_k {x(x-1)\cdots (x-k+1)\over k!}t^k \in\mQ[[x,t]]
$$
where by convention $x(x-1)\cdots (x-k+1)=1$ if $k=0.$ If 
$R$ is a $\mQ$-algebra and ${\mathfrak n}\in R$ is a nilpotent element, then the expression 
$$
\sum_k {x(x-1)\cdots (x-k+1)\over k!}{\mathfrak n}^k
$$ defines a polynomial with coefficients in $R$ and we shall write 
$\TT(x,{\mathfrak n})\in R[x]$ for this polynomial.

We  record the following elementary lemma, which will be used a number of times.

\begin{lemma}
Let $H$ be a vector space over  an infinite field $K$. Let $b_0,\dots, b_l\in H.$ 
Suppose that $\sum_{i=0}^l t^i b_i=0$ for infinitely many $t\in K$. Then 
$b_0=b_1=\dots= b_l=0$.
\label{lemVS}
\end{lemma} 
\beginProof (of Lemma \ref{lemVS}). Let $\{h_c\}_{c\in B}$ be a basis of $H$. Write $b_i=\sum_{c\in B} b_{ic} h_c.$ 
By construction, we have $\sum_{i=0}^l b_{ic} t^i=0$ for all $c\in B$ and infinitely many 
$t\in K$. 
Since $K$ is an infinite field, we thus have $b_{ic}=0$ for all $c\in B$. 
Hence $b_0=b_1=\dots= b_l=0$.
\endProof

For the definition of a family of ample line bundles, which is used in the next lemma, see \cite[I 2.2.3]{SGA6}. We recall that if a scheme is noetherian, separated, and quasi-projective over an affine scheme, then it carries a family of ample line bundles. Also, a separated, noetherian and regular scheme carries an ample family of line bundles.

If $V$ is a vector bundle, we shall write $\rk(V)$ for its rank. 

\begin{lemma} 
\begin{itemize}
\item[{\rm (a)}] Suppose that $R$ is a $\mQ$-algebra and that $r-1\in R$ is nilpotent. Then the identity  $$r^{\otimes n}=\TT(n,r-1)$$
holds in $R$ for all $n\geq 0$ and we have 
$$
\TT(x,r-1)=\sum_{j=0}^{\delta}\big[\sum_{k=0}^{\delta}\sum_{u=0}^k\ {(-1)^{k-u}s(k,j)\over u!(k-u)!}r^u\big] x^j
$$
in $R$ for any $\delta\geq 0$ such that $(r-1)^{\delta+1}=0$ in 
$R$. 

\item[\rm (b)] Let $R$ be a commutative ring and let $r\in R$. Let $r_0\geq 1$. 
Suppose that $r-r_0$ is nilpotent. Then $r$ is invertible in $R[{1\over r_0}]$. 

\item[{\rm (c)}] Let $W$ be a scheme. Let $V$  be a vector bundle  on $W$. 
\begin{itemize}\item[{\rm (1)}] Suppose that $W$ is noetherian. Then $\rk(V)-V$ is nilpotent in $\r K^0(W)$.

\item[{\rm (2)}] Suppose that $W$ is noetherian and carries an ample family of line bundles. Then $\rk(V)-V$ is nilpotent in $K^0(W)$. In particular, there is an isomorphism of rings $K^0(W)/{\rm nilradical}(K^0(W))\simeq \mZ$.
 
\item[{\rm (3)}] If $W$ is noetherian, has an ample line bundle and 
$\dim(W)<\infty$  then we  have 
$(\rk(V)-V)^{\otimes (\dim(W)+1)}=0$ in $K^0(W)$. 
\end{itemize}
\end{itemize}

\label{lemPOL}
\end{lemma}

\beginProof (a)  We compute 
\begin{eqnarray*}
&&r^n=(1+(r-1))^n=
\sum_{k=0}^n{n\choose k}(r-1)^k=\sum_k {n(n-1)\cdots (n-k+1)\over k!}(r-1)^k\\&=&\sum_{k=0}^\delta {n(n-1)\cdots (n-k+1)\over k!}(r-1)^k=\TT(n,r-1)
\end{eqnarray*}
and from the definition of the Stirling numbers of the first kind, we have 
\begin{eqnarray*}
&&\TT(x,r-1)=\sum_{j=0}^{\delta}\big[\sum_{k=0}^{\delta}{1\over k!}s(k,j) (r-1)^k]x^j=
\sum_{j=0}^{\delta}\big[\sum_{k=0}^{\delta}\sum_{u=0}^k\ {(-1)^{k-u}s(k,j)\over u!(k-u)!}r^u\big] x^j.
\end{eqnarray*}
which establishes (a).

(b) We have the identities
$$
{1\over x}={1\over r_0-(r_0-x)}={r_0^{-1}\over 1-{r_0-x\over r_0}}=
r_0^{-1}\sum_k {(r_0-x)^k\over r_0^k}
$$
and setting $x=r$ we get (b).

(c) (1) By noetherian induction. We may thus assume that the claim holds for any proper 
closed subscheme of $W$ instead of $W$. Let $j:U\hookrightarrow W$ be an open subset 
such that $V|_U\simeq\CO_U^{\oplus\rk(V)}$. Let $i:Z\hookrightarrow W$ be the complement of $U$ (viewed as a reduced closed subscheme). If 
$Z=\emptyset$, there is nothing to prove so we may assume that $Z\not=\emptyset.$ 
Recall that there is an exact sequence 
$$
K_0(Z)\stackrel{i_*\,\,\,\,}{\to} K_0(W)\stackrel{j^*}{\to} K_0(U)\to 0
$$
(see \cite[Prop. 3.2]{Quillen-Higher}). Here $i_*$ sends a coherent sheaf on $Z$ to its direct image on $W$ and $j^*$ sends a coherent sheaf on $W$ to its restriction to $U$. 

Let $y\in K_0(W)$. By construction, we have
$
j^*((V-\rk(V))\otimes y)=0
$ and hence 
$$(V-\rk(V))\otimes y=i_*(y_1)$$ for some $y_1\in K_0(Z).$ So if $n>0$ we  have 
\begin{equation}
(V-\rk(V))^{n}\otimes y=i_*((V|_Z-\rk(V))^{\otimes(n-1)}\otimes y_1)
\label{eqINTE}
\end{equation} by the projection formula.  By the inductive hypothesis, 
there is an $n_0\geq 1$, which is independent of $y_1$ and which is such that $(V|_Z-\rk(V))^{\otimes n_0}\otimes y_1=0.$ We conclude from \refeq{eqINTE} that $(V-\rk(V))^{n_0+1}\otimes y=0.$ 
Since $y\in K_0(W)$ was arbitrary, we see that 
$V-\rk(V)$ is nilpotent in $\r K^0(W).$

(c) (2) For this, see \cite[Lemme 1.6, proof]{Thom-L}. 

(c) (3) For this, see \cite[VI, Prop. 6.1]{SGA6}.\endProof

We start with the following proposition, from which equality \refeq{eqRDE} will be deduced.

If $H$ is a complex vector space, we shall say that a topology on $H$ is a vector space topology if it is induced by a norm (note that every complex vector space can be endowed with a norm, so this always exists).

\begin{prop}
Let $T$ be a $\mC$-algebra such that $T\not=0$. Suppose that the natural homomorphism of $\mC$-algebras 
$\mC\to T/{\rm nilradical(T)}$ is surjective (and hence bijective). 

Let \mbox{$\phi(t)=\sum_k  h_k t^k\in T[[t]]$} 
{\rm (} resp.  $\psi(t)=\sum_k  g_k t^k\in T[[t]]$\,{\rm )}.  
Suppose that there is a polynomial $$P(x)=\sum_j d_j x^j\in T[x]$$
{\rm (}\, resp. $$Q(x)=\sum_j e_j x^j\in T[x]\,\,{\rm )}$$ such that $h_k=P(k)$ (resp. $g_k=Q(k))$. 
Then

{\rm (a)} For any $y\in\mC$ such that $|y|<1$ and any vector space topology on $T$, the 
series $\phi(y):=\sum_k h_k y^k$ {\rm(} resp.  $\psi(y):=\sum_k g_k y^k$\,{\rm )} 
 {\rm(} resp.  $(\phi\cdot\psi)(y):=\sum_k\,(\sum_{s+r=k} h_s g_r) y^k$\,{\rm )} converges and we have 
 \begin{equation}
(\phi\cdot\psi)(y)=\phi(y)\cdot\psi(y).
\label{eqSEQ}
\end{equation}
and
\begin{equation}
\lim_{y\to -1^+}\phi(y)={1\over 2}\sum_j E_j d_j.
\label{eqEJ}
\end{equation}
{\rm (b)} Let $X$ be a noetherian scheme with an ample family of line bundles.  Let $V$ be a vector bundle on $X$. 
\begin{itemize} 
\item[\rm (1)] There is a  polynomial $P_V(x)\in K^0(X)_\mQ[x]$ such that 
$P_V(k)=\Sym^k(V)$  in $K^0(X)_\mQ$ for all $k\geq 0$. 
\item[\rm (2)] For any vector space topology on $K^0(X)_\mC$ we have
$$
\lim_{y\to -1^+}\sum_k \Sym^k(V) y^k=\Theta^2(V)^{-1}.
$$
in $K^0(X)_\mC.$ 
\end{itemize}
\label{propKEYI}
\end{prop}
Note that by Lemma \ref{lemVS}, the polynomials $P(x)$ and $Q(x)$ are the unique polynomials such that 
$h_k=P(k)$ and $g_k=Q(k)$ for all $k\geq 0.$
\beginProof Write $\rho:T\to\mC$ for the homomorphism of $\mC$-algebras, which is the 
quotient map $T\to T/{\rm nilradical(T)}$ composed with the inverse of the isomorphism 
of $\mC$-algebras  $\mC\to T/{\rm nilradical(T)}$.

(a) Let $T_0\subseteq T$ be the $\mC$-subalgebra 
generated by all the $d_j$ and all the $e_j$. We contend that $T_0$ is a finite dimensional vector space. For this, note that for each $j$, there exists by construction $r_j\geq 1$ such that $(d_j-\rho(d_j))^{r_j}=(e_j-\rho(g_j))^{r_j}=0$. Hence $T_0$ is spanned as $\mC$-vector space by a finite number of monomials in the $d_j$ and $e_j$. 

Now equip $T_0$ with a vector space norm, which makes $T_0$ into a Banach algebra. This can be achieved as follows. 
Consider $T_0$ as a subalgebra of its algebra of linear operators $\End_\mC(T_0)$ by sending an element $e\in T_0$ 
to the operator $(\cdot)\otimes e$. If we choose an arbitrary vector space norm on $T_0$, we have a corresponding 
operator norm on $\End(T_0)$. The norm on $T_0$, which is inherited from such an operator norm then makes 
$T_0$ into a Banach algebra. Note that the topology of $T_0$ does not depend on the norm, since $T_0$ is finite-dimensional.

Let ${\mathfrak T}_0\subseteq T_0[[t]]$ be the vector subspace of power series $\sum_k c_k t^k$ such that 
the sum $\sum_k c_k y^k$ converges absolutely in $T_0$ (for the Banach algebra norm) for all $y\in\mC$ in the open unit disk. 
A straightforward generalisation of Mertens's theorem for Cauchy products of absolutely convergent series  implies that the subspace ${\mathfrak T}_0$ is a subring of $T_0[[t]]$. Mertens's theorem also implies the following. For any $y\in\mC$ in the open unit disk, let $$
\Ev_y:{\mathfrak T}_0\to T_0
$$
be the map defined by the formula $\Ev_y(\sum_k c_k t^k)=\sum_k c_k y^k$. Then $\Ev_y$ is a homomorphism of $\mC$-algebras. 

We will now make use of the following facts from complex analysis. 

For any $j\geq 0$ we have an identity of power series
$$
[t{{\rm d}\over {\rm d}t}]^{\otimes j}({t\over 1-t})=
\sum_k k^j t^k
$$
where $[t{{\rm d}\over {\rm d}t}]^{\otimes j}$ is the operator $t{{\rm d}\over {\rm d}t}$ composed $j$-times with itself (this can be proved by induction on $j$).  Here by convention we have $k^j=0$ if $k=j=0$. Since the Taylor series of 
\mbox{$t/(1-t)$} has radius convergence $1$ around $0$, we conclude from this that $\phi(t),\psi(t)\in{\mathfrak T}$.
In view of the multiplicativity of $\Ev_y(\cdot)$, this already proves equation \refeq{eqSEQ}. 

Now we also have
$$
([t{{\rm d}\over {\rm d}t}]^{\otimes j}({t\over 1-t}))|_{t=-1}=\Li_{-1}(-j)
$$
where
$$
\Li_{-1}(z)=\sum_{k\geq 1}{(-1)^k\over z^k}
$$
is an instance of a polylogarithmic function and also of a Lerch z\^eta function (see 
\cite[chap. XIII]{WW-CMA}). In particular, we have 
\begin{equation}
\Li_{-1}(-j)=\zeta_\mQ(-j)(2^{1+j}-1)={(-1)^{j+1}\over 2}E_j
\label{eqLI}
\end{equation}
and so $\Li_{-1}(-j)={E_j\over 2}$ if $j>0$ 
(recall that $E_j=0$ if $j$ is odd and $>0$). 

Hence we may compute
\begin{eqnarray*}
&&\Ev_y(\phi(t))=\sum_k h_k y^k=\sum_k(\sum_j d_j k^j) y^k=
d_0+\sum_j d_j [t{{\rm d}\over {\rm d}t}]^{\otimes j}({t\over 1-t})|_{t=y}
\end{eqnarray*}
for any $y\in\mC$ in the open unit disk. Using \refeq{eqLI}, we may then calculate
\begin{eqnarray*}
&&\lim_{y\to -1^+}\Ev_y(\phi(t))={1\over 2}\sum_j E_j d_j
\end{eqnarray*}
establishing equation \refeq{eqEJ}. 

We now establish (b) (1). Let $\pi:P:=\Proj(\oplus_k\Sym^k(V))\to X$. Let 
$\CO(1)$ be the canonical line bundle on $P$ (corresponding to the trivial line bundle 
of weight one on the graded sheaf of algebras $\oplus_k\Sym^k(V)$). Recall that we have $\Sym^k(V)=\pi_*(\CO(k))=\R^\bullet \pi_*(\CO(k))$. Note also that $P$ also carries a family of ample 
line bundles, because $\CO(1)$ is relatively ample (see \cite[2.12 (f)]{TT}). Hence, by Lemma \ref{lemPOL} (a) and (c) (2),  
there is a polynomial $P_{\CO(1)}(x)=\sum_j a_j x^j\in K^0(P)_\mQ[x]$ such that for all $n\geq 0$ we have 
 $P_{\CO(1)}(n)=\CO(n)$ in $K^0(P)_\mQ$. 
 Hence we have 
 $$
 \R^\bullet\pi_*(\CO(n))=\sum_j \R^\bullet\pi_*(a_j )n^j=\Sym^n(V)
 $$
 in $K^0(X)_\mQ$.
 
We now turn to (b) (2). Note first that by Lemma \ref{lemPOL} (c) (2), the nilradical 
of $K^0(X)$ is the kernel of the rank function. In particular, we have 
an exact sequence
$$
0\to {\rm nilradical}(K^0(X))_\mC\to K^0(X)_\mC\to \mC\to 0.
$$
A simple calculation shows that ${\rm nilradical}(K^0(X))_\mC\subseteq 
{\rm nilradical}(K^0(X)_\mC)$. Since $\mC$ is a field, this implies that ${\rm nilradical}(K^0(X))_\mC=
{\rm nilradical}(K^0(X)_\mC)$ and so we may apply (a) and set $T:=K^0(X)_\mC.$ 
Let now $\Sym_t(V):=\sum_k \Sym^k(V)t^k$ and $\Lambda_t(V):=\sum_k \Lambda^k(V)t^k$. Note that $\Lambda_t(V)$ is polynomial in $t$. By standard properties of symmetric and exterior powers, we have
$$
\Sym_t(V)\Lambda_{-t}(V)=1
$$
and in particular
$$
\Sym_t(V)=(\Lambda_{-t}(V))^{-1}
$$
in $K^0(X)_\mQ[[t]].$ Applying (a), we compute \begin{eqnarray*}
&&\lim_{y\to -1^+}(\Sym_y(V)\Lambda_{-y}(V)))=(\lim_{y\to -1^+}\Sym_y(V))
(\lim_{y\to -1^+}\Lambda_{-y}(V))\\&=&(\lim_{y\to -1^+}\Sym_y(V))(1+V+\Lambda^2(V)+\dots+
\Lambda^{\rk(V)}(V))=1
\end{eqnarray*}
in $T=K^0(X)_\mC$, which establishes (b) (2).\endProof

\begin{rem}\rm In the situation of Proposition \ref{propKEYI} (b), suppose that $V$ is a line bundle and that 
$X$ is noetherian and has an ample line bundle.  There is then a Chern character homomorphism $\ch:K^0(X)_\mQ\to{\rm Gr}K^0(X)_\mQ$, 
where ${\rm Gr}K^0(X)_\mQ$ is the graded ring associated with the $\gamma$-filtration on $K^0(X)$ (we refer to \cite[Exp. V \&  VIII]{SGA6} for details). The map $\ch(\cdot)$ is an isomorphism and has the same formal properties as 
the Chern character with values in the Chow intersection ring (which is only defined 
under more restrictive conditions on $X$). Now let $P_V(x)=:\sum_j a_j x^j\in K^0(X)_\mQ$. We may compute
$$
\ch(V^{\otimes k})=\ch(\Sym^k(V))=\sum_j {k^j\over j!}\ch_1(V)^j
$$
In particular, we have $\ch(a_j)={1\over j!}\ch_1(V)^j$. Hence 
$$
\ch((1+V)^{-1})=(1+e^{\ch_1(V)})^{-1}={1\over 2}\sum_j {E_j\over j!}\ch_1(V)^j={1\over 2}\sum_j E_j\ch(a_j)
$$
and since the Chern character $K^0(X)_\mQ\to{\rm Gr}K^0(X)_\mQ$ is an isomorphism, we get 
$$
(1+V)^{-1}={1\over 2}\sum_j E_j a_j
$$
which in combination with (a) proves the formula in (b) (2) (under the just described supplementary assumptions). It seems difficult to generalise this method proof to the case where 
$\rk(V)>1$ however, because there are no explicit formulae for the Chern character of symmetric powers.\end{rem}

The next corollary contains in particular a proof of equality \refeq{eqRDE}.

\begin{cor}
Let $X$ be a scheme which carries an ample family of line bundles.
Let $V$ be a vector bundle. Let $\pi:P:=\Proj(\bigoplus_k\Sym^k(V))\to X$. Let 
$\CO(1)$ be the canonical line bundle on $P$ (corresponding to the trivial line bundle 
of weight one on the graded sheaf of algebras $\bigoplus_k\Sym^k(V)$). Then we have the equality 
\begin{equation}
\R^\bullet\pi_*((1+\CO(1))^{-1})=\Theta^2(V)^{-1}
\label{eqO1}
\end{equation}
in $K^0(X)_\mQ.$ Furthermore, we have 
\begin{equation}
\Theta^2(V)^{-1}={1\over 2}\sum_{j=0}^\delta E_j \sum_{k=0}^{\delta}\sum_{u=0}^k\ {(-1)^{k-u}s(k,j)\over u!(k-u)!}\Sym^u(V)
\label{eqSYMEQ}
\end{equation}
in $K^0(X)_\mQ$ for any $\delta$ such that $(\CO(1)-1)^{\otimes(\delta+1)}=0$ in $K^0(P)_\mQ$. 
\label{lemMIR}
\end{cor}
\beginProof The argument for proving \refeq{eqO1} already appears in the proof of Proposition \ref{propKEYI} (b). Let $$P_{\CO(1)}(x)=\sum_j a_j x^j\in K^0(P)_\mQ[x]$$ be a  polynomial such that 
$\CO(n)=P_{\CO(1)}(n)$ in $K^0(P)_\mQ$ for all $n\geq 0$.
Then we have 
\begin{equation*}
\Sym^n(V)=\R^\bullet\pi_*(\CO(n))=\sum_j(\R^\bullet\pi_*(d_j)) n^j
\end{equation*}
in $K^0(X)_\mQ$. On the other hand, by Proposition \ref{propKEYI} (a) \& (b), we have
$$
{1\over 2}\sum_j E_j d_j=(1+\CO(1))^{-1}
$$
 in $K^0(P)_\mQ$ and 
\begin{equation}
{1\over 2}\sum_j E_j\R^\bullet\pi_*(d_j)=\Theta^2(V)^{-1}
\label{eqIII}
\end{equation}
in $K^0(X)_\mQ$. This proves equality \refeq{eqO1}. For equality \refeq{eqSYMEQ}, note that by Lemma \ref{lemPOL}, we have 
$$
P_{\CO(1)}(x)=\sum_{j=0}^{\delta}\big[\sum_{k=0}^{\delta}\sum_{u=0}^k\ {(-1)^{k-u}s(k,j)\over u!(k-u)!}\CO(u)\big] x^j.
$$
Hence 
$$
\Theta^2(V)^{-1}={1\over 2}\sum_{j=0}^\delta E_j\sum_{k=0}^{\delta}\sum_{u=0}^k\ {(-1)^{k-u}s(k,j)\over u!(k-u)!}\Sym^u(V)
$$
by \refeq{eqIII}.\endProof

\begin{rem}\rm The first part of Corollary \ref{lemMIR} looks like a consequence of the Grothendieck-Riemann-Roch theorem (for the morphism $\pi$). We were not able to deduce Corollary \ref{lemMIR} from this theorem however. We propose this to the 
reader as a challenge. 
\end{rem}

\section{The geometric fixed point formula for an involution}
\label{secFP}

In this section, we review Thomason's geometric fixed point formula and we compute its local 
term in the situation where the group action is of order $2$. For more details, see 
Thomason's articles \cite{Thom-L} and \cite{Thom-Alg}.

Let $S$ be a scheme. Let 
$X/S$ be a $S$-scheme. Let $G$ be a flat group scheme over $S$. 
 Suppose that $X$ is endowed with a $G$-action over $S$. We shall call such a $S$-scheme 
 an equivariant $S$-scheme. 
We shall write $K_0(X,G)$ (resp. $K^0(X,G)$) for the Grothendieck group of the abelian category of $G$-equivariant 
coherent sheaves (resp. the additive category of $G$-equivariant locally free coherent sheaves) on $X$. 
As usual, the tensor product endows $K^0(X,G)$ with a natural structure of commutative ring. 
There is a natural "forgetful" map of groups $K^0(X,G)\to K_0(X,G)$, which sends a $G$-equivariant locally free coherent sheaf on its class in  $K_0(X,G)$. Furthermore, the tensor product by equivariant coherent locally free sheaves endows $K_0(X,G)$ with a structure of $K^0(X,G)$-module. 

If $G$ is the trivial group scheme, we shall 
drop the reference to $G$ in the above (so $K_0(X,G)$ will be written 
$K_0(X)$ etc.). Note that if $X$ has a trivial $G$-equivariant structure, then 
there is a natural morphism of abelian groups (resp. of rings) $K_0(X)\to 
K_0(X,G)$ (resp. $K^0(X)\to 
K^0(X,G)$) which sends a coherent sheaf (resp. a locally free coherent sheaf) 
to the corresponding sheaf endowed with its trivial equivariant structure. 

Suppose from now on that $2$ is invertible on $S$. 
Let $\mu_2=\mZ[T]/(T^2-1)$ be the group scheme of  square roots of unity. 
We set $G:=\mu_{2,S}$ (so that $G\simeq(\mZ/2\mZ)_S$). 

If the action on $X$ 
is trivial and 
$F$ is an equivariant coherent sheaf on $X$, we shall write 
$F_{-}$ for the weight one part of $F$ and $F_+$ for the weight $0$-part 
(in particular, we have $F=F_{+}\oplus F_{-}$).

If $Y$ is a noetherian 
$G$-equivariant $S$-scheme and $f:X\to Y$ is an equivariant and proper $S$-morphism then there is a natural morphism $\R^\bullet f_*:K_0(X,G)\to K_0(Y,G)$, which sends 
a coherent $G$-equivariant sheaf $F$ to the element 
$
\sum_{k\geq 0}(-1)^k\R^k f_*(F)
$
of $K_0(Y,G)$. Here the functors $\R^k f_*$ are the right derived functors of the functor $f_*$  and the sheaves $\R^k f_*(F)$ are endowed with their natural equivariant structure. 
If $f$ is a finite morphism, we shall often write $f_*$ for $\R^\bullet f_*$.

If $f:X\to Y$ is an equivariant $S$-morphism 
of equivariant $S$-schemes, the operation of pull-back of locally free sheaves induces a natural 
map of commutative rings $f^*:K^0(Y,G)\to K^0(X,G)$. 

If $f$ is of finite tor-dimension, 
there is a map of groups $\L^\bullet f^*:K_0(Y,G)\to K_0(X,G)$ which sends 
a coherent $G$-equivariant sheaf $F$ on $Y$ to the element 
$
\sum_{k\geq 0}(-1)^k\L^k f^*(F)
$
of $K_0(X,G)$. Here the functors $\L^k f^*$ are the left derived functors of the functor $f^*$  and the sheaves $\L^k f^*(F)$ are endowed with their natural equivariant structure. The map 
\mbox{$\L^\bullet f^*:K_0(Y,G)\to K_0(X,G)$} is naturally compatible with the map $f^*:K^0(Y,G)\to K^0(X,G)$ 
via the forgetful map.

If $f$ is proper and $Y$ is noetherian, there is a projection formula
$$
\R f_*(x\otimes f^*(e))=\R f_*(x)\otimes e
$$
for any element $x\in K_0(X,G)$ and any element $e\in K^0(Y,G)$. 

We shall write $X_G\hookrightarrow X$ for the fixed point scheme of $X$ (when it exists), 
which is a closed equivariant subscheme of $X$. The scheme $X_G$ has a trivial equivariant structure and it represents the functor 
on $S$-schemes $T\mapsto X(T)^{G(T)}$, where $X(T)$ is the subset of $X(T)$, which 
is fixed under the action of $G(T)$. In our situation, $X_G$ exists if $X$ is separated over $S$ and it is then given by the intersection between the diagonal of $X\times_S X$ and the 
graph of the automorphism of $X$ induced by $1\in\mZ/2\mZ.$ 

We shall write $R(G):=\mZ[\mZ/2\mZ]\simeq\mZ[x]/(x^2-1)$ for the 
group algebra of $\mZ/2\mZ$. There is a unique morphism of rings 
$R(G)\to\mQ$, which sends $x$ to $-1$, and we shall often 
view $\mQ$ as a $R(G)$-algebra via this morphism. 
The ring $K^0(S,G)$ also has a natural $R(G)$-algebra structure, which 
arises from the morphism of rings $R(G)\to K^0(X,G)$ which sends 
$x$ to the structure sheaf $\CO_S$ endowed with a weight one equivariant structure.
Via the natural pull-back maps $K^0(S,G)\to K^0(X,G)$, this endows 
$K^0(X,G)$ with a natural structure of $R(G)$-algebra and $K_0(X,G)$ with a structure of 
$R(G)$-module. This structure 
is compatible with the push-forward and pull-back maps $\R^\bullet f_*$, $f^*$ and $\L^\bullet f^*$.

If $E$ is an equivariant coherent sheaf on $X$, we shall write $E\{1\}$ for the tensor 
product of $E$ with the trivial sheaf endowed with its equivariant structure of weight $1$. Note that by construction, we have the equality $E\{1\}=-E$ in the group 
$K_0(X,G)\otimes_{R(G)}\mQ$. If the action of $G$ on $X$ is trivial and $E$ is an coherent equivariant sheaf on $X$, we shall write 
$[E]_\triv$ for the coherent sheaf $E$ endowed with its trivial equivariant structure. 

Suppose for the time of the next paragraph that the action of $G$ on $X$ is trivial. 
There is then natural map of $\mQ$-vector spaces $\Tr:K_0(X,G)\otimes_{R(G)}\mQ\to K_0(X)_\mQ$, 
which sends $E\otimes r$ to $(E_{+}-E_{-})\otimes r$. To check this, we only have to check that 
this map is $R(G)$-bilinear. This follows from the equalities 
$$
((E\{1\})_{+}-(E\{1\})_{-})\otimes r=(E_- - E_+)\otimes r=(-(E_{+}-E_{-}))\otimes r=
(E_{+}-E_{-})\otimes (-r).
$$
Note that that the map $\Tr$ splits the natural map $K_0(X)_\mQ\to K_0(X,G)\otimes_{R(G)}\mQ$. 
In particular, the natural map $K_0(X)_\mQ\to K_0(X,G)\otimes_{R(G)}\mQ$ is injective.
Also, from the definition we see that $\Tr(\cdot)$ is a map of $K^0(X)_\mQ$-modules.

\begin{theor} Suppose that $S$ is a separated noetherian scheme. 
Suppose that $X$ is a $G$-equivariant $S$-scheme, which is separated and of finite type over $S$. 

{\rm (a)} The fixed point scheme $\iota:X_G\hookrightarrow X$ exists. Write $N_{X_G/X}$ for the conormal 
bundle of $X_G$ in $X$. If every point of $X_G$ has a $G$-invariant affine open neighbourhood  which lies over an open affine subscheme of $S$, then $N_{X_G/X}=N_{X_G/X,-}$.

{\rm (b) [Thomason]} The map $\iota_*\otimes\Id_{\mQ}:K_0(X_G,G)\otimes_{R(G)}\mQ
\to K_0(X,G)\otimes_{R(G)}\mQ$ is an isomorphism of $R(G)$-modules.

\label{theorThom}
\end{theor}
The (unique) element $e\in K_0(X_G,G)\otimes_{R(G)}\mQ$ such that $(\iota_*\otimes\Id_{\mQ})(e)=1$
 will be called the {\it local term} of the equivariant $S$-scheme $X$. 
 
\beginProof (a) The existence of the fixed point scheme in our situation has already been discussed. For the second statement, note that by assumption, we may assume that $X$ and 
$S$ are affine. So suppose that $X=\Spec(B)$ and $S=\Spec(A).$ The action of $G$ on $X$ is given by a $\mZ/2\mZ$ ring grading 
$B=B_{+}\oplus B_-$. The ideal of $X_G$ is the intersection of 
the diagonal of $X\times_S X$ with the graph of $1\in\mZ/2\mZ$. The ideal of 
the diagonal of $X\times_S X=\Spec(B\otimes_A B)$ is generated by the elements 
$b\otimes 1-1\otimes b$, where $b\in B$. This ideal is the 
kernel of the morphism $\nabla:B\otimes_A B\to B$ of $A$-algebras sending $b\otimes c$ to $bc$ (this morphism corresponds to the diagonal immersion). Thus the ideal of the graph of $1\in\mZ/2\mZ$ is 
generated by the elements 
$$
(b_++b_-)\otimes 1-1\otimes(b_+-b_-)=(b_+\otimes 1-1\otimes b_+)+(b_-\otimes 1+1\otimes b_-).$$ 
The image of this last ideal under $\nabla$ is the set of element $2b_-$, with $b_-\in B_-$. 
Since $2$ is a unit, this is exactly the ideal $(B_-)$ generated by $B_-$. We thus 
see that the ideal of $X_G$ in $X$ is $(B_-)$. In particular, 
$(B_-)/(B_-)^2$ (which corresponds to $N_{X_G/X}$) is generated by homogenous elements of weight one. This proves the claim. 

(b) This is a special case of \cite[Th. 2.2]{Thom-L}. \endProof

\begin{theor}
Suppose that $S$ is a separated noetherian scheme and that $X$ is a $G$-equivariant $S$-scheme, which is separated and of finite type over $S$. Suppose 
that $\iota:X_G\hookrightarrow X$  is a Cartier divisor in $X$. Suppose furthermore that every point of $X_G$ has an open affine $G$-invariant neighbourhood, which lies over an open affine subscheme of $S.$ 

Then $N_{X_G/X}-1\in \r K^0(X)_\mQ$ is nilpotent. Furthermore, if we let 
$$
\sum_j a_j x^j:=\TT(x,N_{X_G/X}-1)\in \r K^0(X_G)_\mQ[x]
$$
Then   we have 
\begin{equation}
(\iota_*\otimes\Id_\mQ)\big[{1\over 2}\sum_j E_j a_j\big]=1
\label{eqTC}
\end{equation}
in $K_0(X,G)\otimes_{R(G)}\mQ$. In other words, the local term of $X$ is ${1\over 2}\sum_j E_j a_j$. \label{theorEulerCartier}
\end{theor}
Note that the element ${1\over 2}\sum_{j} a_jE_j$ lives in $\r K^0(X_G)_\mQ$ 
but we take its image in \mbox{$K_0(X_G,G)\otimes_{R(G)}\mQ$} in the last equality.

\beginProof (of Theorem \ref{theorEulerCartier}). The fact that $N_{X_G/X}-1\in \r K^0(X)_\mQ$ 
is nilpotent is a consequence of Lemma \ref{lemPOL}. 

By Theorem \ref{theorThom} (b), there is a unique element \mbox{$e\in K_0(X_G,G)\otimes_{R(G)}\mQ$} such that 
$(\iota_*\otimes\Id_\mQ)(e)=1$. By the projection formula, we have 
$$
(\iota_*\otimes\Id_\mQ)(e\otimes\CO(-n X_G)|_{X_G})=\CO(-n X_G)
$$
in $K_0(X,G)\otimes_{R(G)}\mQ$. We have
$$
\CO(-n X_G)|_{X_G}=(-1)^n(\CO(-X_G)|_{X_G}\{1\})^{\otimes n}
$$ in $K_0(X_G,G)\otimes_{R(G)}\mQ$ and thus the adjunction formula implies that  we have 
$$
(-1)^n(\iota_*\otimes\Id_\mQ)(e\otimes(\CO(-X_G)|_{X_G}\{1\})^{\otimes n})=(-1)^n(\iota_*\otimes\Id_\mQ)(e\otimes(N_{X_G/X}\{1\})^{\otimes n})=\CO(-n X_G)
$$
in $K_0(X,G)\otimes_{R(G)}\mQ$. Note that by Theorem \ref{theorThom} (a), the sheaf 
$N_{X_G/X}\{1\}$ has a trivial equivariant structure. 

On the other hand, we have an exact equivariant sequence
$$
0\to \CO(-n X_G)\to \CO_X\to \CO_{X_G^{(n-1)}}\to 0
$$
where $X_G^{(n-1)}$ is the $n-1$-th infinitesimal neighbourhood of $X_G$ in $X$. 
Since the immersion of $X_G$ into $X$ is regular, we have an equality 
$$
\CO_{X_G^{(n-1)}}=\sum_{k=0}^{n-1}(\iota\otimes\Id_\mQ)_*(N_{X_G/X}^{\otimes k})=\sum_{k=0}^{n-1}(-1)^k(\iota\otimes\Id_\mQ)_*((N_{X_G/X}\{1\})^{\otimes k})
$$
in  $K_0(X,G)\otimes_{R(G)}\mQ$. Thus we have 
$$
1=(\iota_*\otimes\Id_\mQ)(e)=\sum_{k=0}^{n-1}(-1)^k(\iota_*\otimes\Id_\mQ)((N_{X_G/X}\{1\})^{\otimes k})+
(-1)^n(\iota_*\otimes\Id_\mQ)(e\otimes(N_{X_G/X}\{1\})^{\otimes n})
$$
and thus by Lemma \ref{lemPOL} we have 
\begin{eqnarray*}
&&1=(\iota_*\otimes\Id_\mQ)(e)=\sum_{k=0}^{n-1}(-1)^k(\iota_*\otimes\Id_\mQ)(
\sum_j a_j k^j)+
(-1)^n(\iota_*\otimes\Id_\mQ)(e\otimes\sum_j a_j n^j))\\
&=&
\sum_j(\iota_*\otimes\Id_\mQ)(
a_j\sum_{k=0}^{n-1}(-1)^k k^j)+
(-1)^n(\iota_*\otimes\Id_\mQ)(e\otimes\sum_j a_j n^j))\\
&=& 
(\iota_*\otimes\Id_\mQ)\Big[\sum_{j=0}^\delta
\Big(\sum_{k=0}^{n-1}(-1)^k k^j)+
(-1)^n n ^j e\Big)\otimes a_j)
\end{eqnarray*}

We now use the following formula from \cite[p. 2]{Kim-Alt}. This formula states that 
$$
\sum_{k=0}^{n-1}(-1)^k k^j={(-1)^{n+1}\over 2}\sum_{l=0}^{j-1}{j\choose l} E_l n^{j-l}+{E_j\over 2}
(1+(-1)^{n+1}).
$$
In view of the injectivity of $\iota_*\otimes\Id_\mQ$, this formula implies that
$$
(-1)^ne=
-{1\over 2}\sum_j\Big[\sum_{l=0}^{j-1}{j\choose l} E_l n^{j-l}+{E_j\over 2}
(1+(-1)^{n+1})\Big]a_j+\sum_j n^j (e\otimes a_j)
$$
for all $n\geq 0$ and hence
$$
(-1)^n\Big[e-{1\over 2}\sum_j {E_j a_j}\Big]=
-{1\over 2}\sum_j\Big[\sum_{l=0}^{j-l}{j\choose l} E_l n^{j-l}+{E_j\over 2}\Big]a_j+\sum_j n^j (e\otimes a_j)
$$
Since the left hand side of the last equality takes the same value for infinitely many 
$n$, we deduce from the Lemma \ref{lemVS} that its right hand side vanishes.
Hence 
$$
e={1\over 2}\sum_j {E_j a_j}
$$
whence the Theorem.
\endProof

\begin{rem}\rm We keep the assumptions of Theorem \ref{theorEulerCartier}. 
Suppose in addition that $X_G$ carries an ample family of line bundles. Then it follows from 
Lemma \ref{lemPOL} (a) and Proposition \ref{propKEYI} (a) \& (b) that ${1\over 2}\sum_j {E_j a_j}$ 
is the image in $K_0(X_G,G)\otimes_{R(G)}\mQ$ of $\Theta^2([N_{X_G/X}]_\triv)^{-1}\in K^0(X_G)_\mQ.$ 
\label{remINVREPC}
\end{rem}

\begin{cor}
Suppose that $S$ is a separated noetherian scheme and that $X$ is a $G$-equivariant $S$-scheme, which is separated and of finite type over $S$. Suppose furthermore that every point of $X_G$ has an open affine $G$-invariant neighbourhood, which lies over an open affine subscheme of $S.$ Let 
$I$ be the ideal sheaf of $X_G$ in $X$. Let $\iota:X_G\to X$ 
be the immersion of $X_G$ into $X$. 

Let $\pi:\wt{X}\to X$ be 
the blow-up of $X$ along $I$ and let $\phi:E\to X_G$ be the corresponding 
exceptional divisor. 

Let $\lambda\geq 0$ be such that 

- for all $r\geq \lambda$ and $a>0$ we have 
$\R^a\pi_*(\CO(-E)^{\otimes r})=0$;

- for all $r\geq \lambda$, the 
morphism of sheaves $I^r\to \pi_*(\CO(-E)^{\otimes r})$ is an isomorphism.

Then for all $\delta\geq\delta_0(N_{E/\wt{X}})$  we have 
\begin{eqnarray}
&&(\iota_*\otimes\Id_\mQ)({(-1)^\lambda\over 2}\sum_{j=0}^\delta E_j\sum_{k=0}^{\delta}\sum_{u=0}^k\ {(-1)^{k-u}s(k,j)\over u!(k-u)!}[I^{u+\lambda}/I^{u+1+\lambda}]_\triv\label{eqFPG}\\&+&
\sum_{k=0}^{\lambda-1}(-1)^k[I^k/I^{k+1}]_\triv)
=1\nonumber
\end{eqnarray}
in $K_0(X,G)\otimes_{R(G)}\mQ$. 
\label{corFPG}
\end{cor}
In other words, the local term of $X$ is given by the formula
$${(-1)^\lambda\over 2}\sum_{j=0}^\delta E_j\sum_{k=0}^{\delta}\sum_{u=0}^k\ {(-1)^{k-u}s(k,j)\over u!(k-u)!}[I^{u+\lambda}/I^{u+1+\lambda}]_\triv+
\sum_{k=0}^{\lambda-1}(-1)^k[I^k/I^{k+1}]_\triv.$$

Note  that we have $\R^\bullet\phi(\CO(-E)|_E^{\otimes r})=(I^r/I^{r+1})^r$ 
(in the derived category) for all $r\geq \lambda.$ To see this, apply $\R^\bullet\pi_*(\cdot)$
to the exact sequence
$$
0\to\CO(-(r+1)E)\to \CO(-rE)\to \CO(-E)|_E^{\otimes r}\to 0.
$$

Finally, note that if $X_G$ is regularly embedded in $X$, then we may take $\lambda=0$. 

\beginProof The existence of $\lambda$ follows from \cite[Cohomology of Schemes, Lemmas 14.2 and 14.3]{StacksProject}. By functoriality, the 
action of $G$ extends to $\wt{X}$ and $E$, making $\pi$ and 
$\phi$ into equivariant morphisms. 

We claim that every point of $\wt{X}_G$ has an open affine $G$-invariant neighbourhood, which lies over an open affine subscheme of $S.$ From the hypothesies, we see that to prove this 
we may assume (for the time of the proof of the claim) that $X=\Spec(B)$ and $S=\Spec(A).$ 
We saw in the proof of Theorem \ref{theorThom} (a) that in this situation 
we have $I=(B_{-}).$ Now we have $\wt{X}=\Proj(\oplus_{k\geq 0} I^k)$ and $\wt{X}$ is covered 
by the open affine subschemes, which are the spectra of the rings 
$(\oplus_{k\geq 0}I^k)_f^{(0)}$, where $f\in B_{-}$. Here $(\oplus_{k\geq 0}I^k)_f^{(0)}$ are the homogenous elements of degree $0$ (for the $\mZ$-grading of $\oplus_{k\geq 0}I^k$) of the localisation of the ring $\oplus_{k\geq 0}I^k$ at $f$. Now $\Spec((\oplus_{k\geq 0}I^k)_f^{(0)})$ is 
$G$-invariant, since $f$ is of of weight one. Hence $\wt{X}$ has an open covering 
by affine $G$-invariant subschemes and so this proves the claim.

We also claim that $E=\wt{X}_G$. Indeed, we have 
$E={\rm Proj}(\oplus_{k\geq 0}I^k/I^{k+1})$ and ${\rm Proj}(\oplus_{k\geq 0}I^k/I^{k+1})$ is covered 
locally by the open affine subschemes, which are the spectra of the rings 
$(\oplus_{k\geq 0}I^k/I^{k+1})_f^{(0)}$, where $f\in I/I^2$. Here, as before, the ring  $(\oplus_{k\geq 0}I^k/I^{k+1})_f^{(0)}$ consists of  the homogenous elements of degree $0$ (for the $\mZ$-grading of $\oplus_{k\geq 0}I^k/I^{k+1}$) of the localisation of the ring $\oplus_{k\geq 0}I^k/I^{k+1}$ at $f$. The elements of 
degree $0$ are all sums of elements of the form $g/f^k$, where $g\in I^k/I^{k+1}$. Such 
elements have weight $(-1)^k/(-1)^k=1$ according to Theorem \ref{theorThom} (a). 
Hence $E$ is fixed under the action of $G$. On the other hand, $\wt{X}_G\subseteq E$ 
since $E$ lies over $X_G$, and hence $E=\wt{X}_G.$

We now apply Theorem \ref{theorEulerCartier} to $\wt{X}$. Let $\delta\geq \delta_0(N_{E/\wt{X}})$ and let 
$$
\sum_{j=0}^\delta a_j t^j=\TT(\rho(N_{E/\wt{X}})-1,t)\in\r K^0(E)_\mQ[t].
$$
From Lemma \ref{lemPOL}, we have 
\begin{equation}
a_j=\sum_{k=0}^{\delta}\sum_{u=0}^k\ {(-1)^{k-u}s(k,j)\over u!(k-u)!}N_{E/\wt{X}}^{\otimes u}.
\label{eqAJ}
\end{equation}
Let $e:E\to {X}$ be the immersion of $E$ into ${X}$.

 We deduce from Theorem \ref{theorEulerCartier} and the adjunction formula that 
\begin{equation}
(e_*\otimes\Id_\mQ)\Big[{1\over 2}\sum_{j=0}^\delta E_j(a_j\otimes N_{E/\wt{X}}^{\otimes \lambda})]=\CO(-E)^{\otimes \lambda}
\label{eqINT}
\end{equation}
(remember the all sheaves in this equation carry their natural equivariant structure). On the other hand, we have an equivariant exact sequence
$$
0\to I^\lambda\to\CO_{X}\to \CO_{X_G^{(\lambda-1)}}\to 0
$$
where $X_G^{(\lambda-1)}$ is the infinitesimal neighbourhood of order $\lambda-1$ of $X_G.$ 
Thus 
$$
I^{\lambda}=1-(\iota_*\otimes\Id_\mQ)[\sum_{k=0}^{\lambda-1}(-1)^k[I^k/I^{k+1}]_\triv]
$$
and putting this together with \refeq{eqINT}  and using the definition of $\lambda$ we have 
\begin{eqnarray*}
&&(\iota_*\otimes\Id_\mQ)(\R^\bullet\phi_*\otimes\Id_\mQ)\Big[{(-1)^\lambda\over 2}\sum_{j=0}^\delta  E_j(a_j\otimes [N_{E/\wt{X}}]_\triv^{\otimes \lambda})\Big]\\&+&(\iota_*\otimes\Id_\mQ)[\sum_{k=0}^{\lambda-1}(-1)^k[I^k/I^{k+1}]_\triv]
=1
\end{eqnarray*}
from which \refeq{eqFPG} follows.\endProof

\begin{rem}\rm 
Keep the assumptions of Corollary \ref{corFPG} and assume that $X_G$ has an ample family of 
line bundles and that $X_G$ is regularly embedded in $X$. It then follows from remark \ref{remINVREPC} 
and Corollary \ref{lemMIR} that the local term \refeq{eqFPG} is equal to the image of 
$\Theta^2([N_{X_G/X}]_\triv)^{-1}$ in $K_0(X,G)\otimes_{R(G)}\mQ$. 
\label{remINVREP}
\end{rem}


\section{The local term for smoothly embeddable schemes}

\label{secREG}

We keep the notation of the last section. In the present section, we shall compute the local term in the situation where the 
equivariant scheme under consideration can be equivariantly embedded by a closed immersion of 
finite tor-dimension into a $G$-equivariant $S$-scheme, whose fixed point scheme is regularly embedded. 

Throughout this section, $S$ will be a noetherian and separated scheme (with $2$ invertible on $S$ as before). 

Suppose that $X$ and $Z$ are $G$-equivariant $S$-scheme which is separated and of finite type over $S$. We suppose that the fixed point scheme $Z_G$ of $Z$ is regularly embedded in $Z$. 
We also suppose given a $G$-equivariant closed immersion $\alpha:X\hookrightarrow Z$, which 
is of finite tor-dimension. In other words, the coherent sheaf $\alpha_*(\CO_X)$ has a 
finite resolution by locally free $\CO_Z$-modules locally on $Z$. 

We let $\iota: X_G\to X$ (resp. $\beta:Z_G\hookrightarrow Z$) be the canonical closed immersion. 
We let $\wt{Z}$ be the blow-up of $Z$ along $Z_G$ and let $E$ be the exceptional 
divisor of $\wt{Z}$.

\begin{prop}
Let $\delta\geq\delta_0(N_{E/\wt{Z}})$.  Then the local term of $X$ is 
$$
{1\over 2}\sum_{j=0}^\delta E_j\sum_{k=0}^{\delta}
\sum_{u=0}^k\ {(-1)^{k-u}s(k,j)\over u!(k-u)!}\Big(\sum_{i\geq 0}(-1)^i\rTor^i_{\CO_Z}(\CO_{Z_G},\CO_X)\otimes_{\CO_{Z_G}}[\Sym^u(N_{Z_{G}/Z})]_\triv\Big)
$$
\end{prop}
Here $\rTor^i_{\CO_Z}(\CO_{Z_G},\CO_X)\otimes_{\CO_{Z_G}}[\Sym^u(N_{Z_{G}/Z})]_\triv$ is viewed 
as a $\CO_{X_G}$-module. Note that $\rTor^i_{\CO_Z}(\CO_{Z_G},\CO_X)$ is a $\CO_{X_G}$-module because it 
is annihilated by the ideal sheaves of $Z_G$ and of $X$, and hence by the ideal sheaf of $Z_G\cap X=X_G.$ 
\beginProof Since $\alpha$ is of finite tor-dimension, we have a homomorphism $\L^\bullet\alpha^*:K_0(Z,G)\to 
K_0(X,G).$ Note also that the local term of $Z$ is given by 
$$
{1\over 2}\sum_{j=0}^\delta E_j\sum_{k=0}^{\delta}\sum_{u=0}^k\ {(-1)^{k-u}s(k,j)\over u!(k-u)!}\Sym^u([N_{Z_G/Z}]_\triv)
$$
by Corollary \ref{corFPG}. 
We now compute 
\begin{eqnarray*}
&&(\L^\bullet\alpha^*\otimes\Id_\mQ)(1)=1\\&=&{1\over 2}\sum_{j=0}^\delta E_j\sum_{k=0}^{\delta}\sum_{u=0}^k\ {(-1)^{k-u}s(k,j)\over u!(k-u)!}(\L^\bullet\alpha^*\otimes\Id_\mQ)((\beta_*\otimes\Id_\mQ)([\Sym^u(N_{Z_{G}/Z})]_\triv))\\&=&
{1\over 2}\sum_{j=0}^\delta E_j\sum_{k=0}^{\delta}
\sum_{u=0}^k\ {(-1)^{k-u}s(k,j)\over u!(k-u)!}\Big(\sum_{i\geq 0}(-1)^i\rTor^i_{\CO_Z}([\Sym^u(N_{Z_{G}/Z})]_\triv,\CO_X)\Big)
\end{eqnarray*}
On the other hand, by \cite[VII, Lemme 2.4]{SGA6}, we have a functorial isomorphism 
$$\rTor^i_{\CO_Z}([\Sym^u(N_{Z_{G}/Z})]_\triv,\CO_X)\simeq\rTor^i_{\CO_Z}(\CO_{Z_G},\CO_X)\otimes_{\CO_{Z_G}}[\Sym^u(N_{Z_{G}/Z})]_\triv$$
and the formula for the local term follows.\endProof

\begin{cor}[{see \cite[Th. 3.5]{Thom-L}} and \cite{Quart-Loc}]
Suppose in addition that $X_G$ carries an ample family of line bundles. 
Then the local term of $X$ is equal to 
$$
\big(\sum_{i\geq 0}(-1)^i\rTor^i_{\CO_Z}(\CO_{Z_G},\CO_X)\big)\otimes_{X_{G}}\Theta^2([N_{Z_G/Z}]_\triv|_{X_G})^{-1}
$$
\label{corGOODLT}
\end{cor}
Here the inverse of $\Theta^2([N_{Z_G/Z}]_\triv)$ is taken in the ring $K^0(X_G)_\mQ.$
\beginProof This follows from the proposition and Corollary \ref{lemMIR}.\endProof

\section{The Castelnuovo-Mumford regularity of Rees algebras}

\label{secCM}

We will now explain the link between the integer $\lambda$ considered in Corollary \ref{corFPG} and the Castelnuovo-Mumford 
regularity of the Rees algebra of the sheaf of ideals of the fixed point scheme. 
The results of this section will be used to show that we can take $\lambda=0$ is certain circumstances 
and can be skipped by a reader not interested in estimating $\lambda.$ 

We need some preparation. Let $T=\oplus_{k\geq 0} T_k$ be a graded 
ring. We view $T$ as a $\mZ$-graded ring. For any $a\geq 0$, we shall write 
$T_{>a}:=\oplus_{k>a} T_k$ (note that this is a graded ideal). Let $M=\oplus_{k\in\mZ}M_k$ be a graded 
$T$-module. For any $i$, we shall write 
$$
H^i_{T_{>0}}(M):=\indlim_{n}{\rm Ext}^i(T/(T_{>0})^n,M)
$$
for the $i$-th local cohomology group of $M$ with respect to $T_{>0}$. This is a $T$-module and 
it inherits a natural $\mZ$-grading from $T$ and $M$ by functoriality. For any $j\in\mZ$, we  shall write $H^i_{T_{>0}}(M)_j$ for the 
set of elements of  $H^i_{T_{>0}}(M)$, which are homogenous of degree $j$. 

Assume now that $M$ is finitely generated and that $T_1$ generates $T$ as a $T_0$-algebra. 
Suppose also that $T_0$ is noetherian. 

The Castelnuovo-Mumford regularity (or simply regularity) of $M$ is then the integer
$$
\reg(M):={\rm max}\{i+j\,|\,H^i_{T_{>0}}(M)_j\not=0\}.
$$
It can be shown (see \cite[Th. 8]{Ooishi-CR}) that $\reg(M)\geq 0$. 

We also recall the following geometrical interpretation of local cohomology in this situation.

Let $\wt{M}$ be the coherent sheaf on $\Proj(T)$ associated with $M$. 
We then have a canonical isomorphism of $T_0$-modules
$$
H^{i+1}_{T_{>0}}(M)_j\simeq H^{i}(\Proj(T),\wt{M}(j))
$$
for all $i\geq 1$ and  $j\in\mZ$, and an exact sequence of $T_0$-modules
$$
0\to H^0_{T_{>0}}(M)_j\to M_j\to H^{0}(\Proj(T),\wt{M}(j))\to H^1_{T_{>0}}(M)_j\to 0
$$
for all $j\in\mZ$ (see \cite[20.4]{Brod-LC}). 

Now consider the following special case. Let $R$ be a noetherian ring 
and let $J\subseteq R$ be an ideal. Let 
$\Rees(J):=\oplus_{k\geq 0} J^k$ be the Rees algebra of $J$, which we view as a $\mZ$-graded 
ring. Let 
$$
H(J):=\inf\{h\in\mZ\,|\,H^{i}_{\Rees(J)_{>0}}(\Rees(J))_q=0\,\,\,\forall i\geq 0\,{\rm and}\,\forall q\geq h\}
$$
Clearly, we have $H(J)\leq\reg(\Rees(J))+1$ and if $R$ is a domain (so that $\Rees(J)$ is also a domain) then we even have $H(J)\leq\reg(\Rees(J))$ (because then we have \mbox{$H^{0}_{\Rees(J)_{>0}}(\Rees(J))=0$}). 

The above definitions tie in with Corollary \ref{corFPG} in the following way.

\begin{lemma}
Assumptions as in Corollary \ref{corFPG}. Suppose that 
$h\geq \sup_{x\in X_G}H(\Rees(I_{\CO_x})).$ Then 

- for all $r\geq h$ and $a>0$ we have 
$\R^a\pi_*(\CO(-E)^{\otimes r})=0$;

- for all $r\geq h$, the 
morphism of sheaves $I^r\to \pi_*(\CO(-E)^{\otimes r})$ is an isomorphism.
\end{lemma}
In particular, any $\lambda\geq \sup_{x\in X_G}H(\Rees(I_{\CO_x}))$ satisfies the condition 
required by Corollary \ref{corFPG}.
\beginProof The proof immediately follows from the definition of the blow-up  and from 
the discussion above.\endProof

We shall now collect a few results from commutative algebra.

We recall the notion of $d$-sequence (see \cite{Huneke-d}). Ideals 
generated by $d$-sequences will turn out to have Rees algebras with vanishing 
Castelnuovo-Mumford regularity in certain circumstances. See below.

Let $R$ be a commutative ring. A sequence $x_1,\dots, x_n\in R$ is said to be 
a $d$-sequence if 

(a) for all $i=1,\dots,n$ we have $x_i\not\in (x_1,\dots, x_{i-1},x_{i+1},\dots, x_n)$;

(b) if $i=1,\dots,n$ and $k\in\{i,\dots,n\}$, then the image of $x_i$ in 
$R/((x_0,\dots,x_{i-1}):(x_k))$ is not a zero divisor.

Here we set $x_0:=0$. Condition (b) can also be written as
$$
((x_0,\dots,x_{i-1}):(x_k x_i))=((x_0,\dots,x_{i-1}):(x_k)).
$$
Note that if $x_1,\dots,x_n$ is a $R$-regular sequence, then it is a $d$-sequence in $R$. 

Write $I:=(x_1,\dots,x_n).$ 
We say that the sequence $x_1,\dots, x_n\in R$ is a Cohen-Macaulay $d$-sequence (see \cite[Def. before Th. 2.3]{Huneke-Koszul}) 
if the rings \mbox{$R/((x_0,\dots,x_{i}):I)$} and \mbox{$R/(((x_0,\dots,x_{i}):I)+I)$} are Cohen-Macaulay rings 
for all $i=0,\dots,n-1$. 

The following lemma (which is probably well-known, but for which we could find no proof in the literature) is needed in the proof of the subsequent proposition.

\begin{lemma}
Let $T$ be a Cohen-Macaulay local noetherian ring with maximal ideal $\Fm_T$. Let $J\subseteq T$ be an ideal. Let $j:=\grade(J)$ and let 
$c=\dim_{T/\Fm_T}(J/\Fm_T J).$ Then $j\leq c$ and there exists 
a set of generators $t_1,\dots,t_c$ of $J$ such that 
$t_1,\dots,t_j$ is a $T$-regular sequence.
\label{lemRSGEN}
\end{lemma}
Recall that $\grade(J)=\grade_T(J)=\depth_T(J,T)$ is the length of a regular $T$-sequence of maximal length contained in $J$.
\beginProof\footnote{I am grateful to the anonymous user 'metalspringpro' on the website math stackoverflow for a comment which suggested this proof.} By induction on $\grade(J)$. If $\grade(J)=0$ there is nothing to prove. 
So suppose that $\grade(J)>0$. Let $\{\Fp_i\}$ be the set of associated primes of the 
the $0$-ideal in $T$. We have $J\not\subseteq\Fp_i$ for all $i$ since $J$ contains a non zero divisor by assumption. We also have $J\not\subseteq\Fm_T J$, for otherwise $J=0$ by Nakayama's lemma. Hence, by prime avoidance, we have 
$$
J\not\subseteq(\cup_i\Fp_i)\cup\Fm_T J.
$$
Hence there exists $t_1\in J$, such that $t_1\not\in (\cup_i\Fp_i)\cup\Fm_T J$.
In particular $t_1$ is not a zero divisor and $t_1\,({\rm mod}\,\Fm_T J)\not=0$. 
Now consider the ideal $J':=J/(t_1)$ in $T':=T/(t_1).$ We have 
$\grade_{T'}(J')<\grade_T(J)$, for otherwise the element $t_1$ together with a lifting to $T$ of a regular sequence of maximal length in 
$J'$ would give a regular sequence of length $>\grade(J)$ in $J$. Thus, by induction, 
there is a minimal set of generators $t_2',\dots,t'_{c}$ of $J'$ such that 
$t_2',\dots,t_{j'+1}$ is a regular $T/(t_1)$-sequence, where $j'=\grade_{T'}(J').$ 
Let $t_2,\dots, t_c$ be liftings of $t_2',\dots,t'_{c}$ to $J$. The sequence $t_1,\dots, t_c$ 
is by construction a minimal set of generators of $J$ and the sequence $t_1,\dots, t_{j'+1}$ is a regular $T$-sequence. 
Furthermore, 
$\height((t_1,\dots, t_{j'+1}))\leq j'+1$ by Krull's principal ideal theorem and 
since $t_1,\dots, t_{j'+1}$ is a $T$-regular sequence, we have $\height((t_1,\dots, t_{j'+1}))\geq j'+1$ (see \cite[before Th.16.9]{Mat-Co}) and thus \mbox{$\height((t_1,\dots, t_{j'+1}))=j'+1$.} Now  by construction $J$ is contained 
in the union of associated prime ideals of $(t_1,\dots, t_{j'+1})$. Also, since $T$ is Cohen-Macaulay, \cite[Cor. 18.14]{Eis-Co} implies that the associated 
prime ideals of $(t_1,\dots, t_{j'+1})$ are all minimal prime ideals of $(t_1,\dots, t_{j'+1})$ and 
thus by prime avoidance $J$ is contained in a minimal prime ideal of $(t_1,\dots, t_{j'+1})$. 
Hence $\height(J)=j'+1$ and by \cite[Th. 17.4]{Mat-Co} we have 
$\grade(J)=j=j'+1$. Hence $t_1,\dots,t_c$ is the required minimal set of generators of $J$.\endProof

\begin{prop}[Huneke]
Let $T$ be a noetherian  local Cohen-Macaulay ring with maximal ideal $\Fm_T$ and let $\Fp\in\Spec(T)$ be a prime ideal of 
$T$. Suppose that $$\dim_{T/\Fm_T}(\Fp/\Fm_T\Fp)\leq\height(\Fp)+1.$$ Suppose also 
that the ring $T_\Fp$ is regular. Then $\Fp$ is generated by a $d$-sequence.
\label{propH0}
\end{prop}
This is \cite[Intro. (5)]{Huneke-Sym} and it is also a special case of \cite[Prop. 2.4 (1)]{Huneke-Koszul}, although neither references provide an explicit proof. We include a proof for the sake of completeness and to show $d$-sequences at work. 

\beginProof (of Proposition \ref{propH0}). 
 Let $n:=\height(\Fp).$ Note that we have $n=\grade(\Fp)$ because 
$T$ is Cohen-Macaulay (see \cite[Th. 17.4]{Mat-Co} ). Let $c:=\dim_{T/\Fm_T}(\Fp/\Fm_T\Fp).$ 
Let $t_1,\dots, t_c$ be a minimal set of generators of $\Fp$ such that 
$t_1,\dots, t_n$ is a $T$-regular sequence of maximal length in $\Fp$ (this exists by the previous lemma). By assumption, 
we have $n\geq c-1$. If $n=c$ then $\Fp$ is generated by a regular $T$-sequence, 
and hence by a $d$-sequence, so there is nothing to prove. So we suppose that $n=c-1$. 
We will show that 
$$
((t_1,\dots,t_{n}):(t_c))=((t_1,\dots,t_{n}):(t_c^2))
$$
thus showing that $t_1,\dots, t_c$ is a $d$-sequence and proving the proposition. 

Note first that $\Fp$ is a minimal prime of $(t_1,\dots,t_{n})$ 
because $\height((t_1,\dots,t_{n}))=n$ (for this last equality see \cite[Th. 17.4]{Mat-Co}). 
Note also that the associated primes of $(t_1,\dots,t_{n})$ are all minimal because 
$T$ is Cohen-Macaulay (for this see \cite[Cor. 18.4]{Eis-Co}). 

Now let $\Fq_1,\dots,\Fq_l$ be a minimal primary decomposition of $(t_1,\dots,t_{n})$. 
We may suppose that $\rad(\Fq_1)=\Fp$ by the above. We have $t_c\not\in\rad(\Fq_i)$ for $i>1$ (and hence $t_c^2\not\in\rad(\Fq_i)$) for otherwise $\Fp\subseteq\rad(\Fq_i)$ and then $\Fp=\rad(\Fq_1)=\rad(\Fq_i)$ (because the associated primes of $(t_1,\dots,t_{n})$ are all minimal), which would contradict the minimality of the decomposition. 

We also claim that $\Fq_1=\Fp.$ To see this, note that since $\Fq_1$ is $\Fp$-primary and 
$T$ is noetherian, the ring $T/\Fq_1$ has only one associated prime ideal, 
namely the ideal $\Fp/\Fq_1$. Also, note that the image of $t_1,\dots,t_n$ in $T_\Fp$ is a regular system of parameters of 
$T_\Fp$ since $T_\Fp$ is regular and of dimension $n$. In particular, $\Fq_{1,\Fp}=(t_1,\dots,t_n)_\Fp$ is the maximal ideal $\Fp_\Fp$ of $T_\Fp$. Now  the localisation of $T/\Fq_1$ at $\Fp/\Fq_1$ is the ring $T_\Fp/\Fq_{1,\Fp}$, and by the above this is a field and in particular  a regular ring.  We conclude from 
\cite[Exercise 11.10, p. 270]{Eis-Co} that $T/\Fq_1$ is reduced, in other words that 
$\Fq_1=\rad(\Fq_1)=\Fp$. 

We may now compute
$$
((t_1,\dots,t_{n}):(t_c))=\cap_i(\Fq_i:(t_c))=(\Fq_1:(t_c))\cap(\cap_{i>1}\Fq_i)=
(\Fp:(t_c))\cap(\cap_{i>1}\Fq_i)=T\cap(\cap_{i>1}\Fq_i)
$$
and 
$$
((t_1,\dots,t_{n}):(t_c^2))=\cap_i(\Fq_i:(t_c^2))=(\Fq_1:(t_c^2))\cap(\cap_{i>1}\Fq_i)=
(\Fp:(t_c^2))\cap(\cap_{i>1}\Fq_i)=T\cap(\cap_{i>1}\Fq_i)
$$
proving the proposition.\endProof

We now come to the relationship between $d$-sequences and the Rees algebra. 

\begin{theor}
Suppose that $T$ is a  local  noetherian ring with infinite residue field and maximal ideal 
$\Fm_T$. Let 
$I\subseteq T$ be an ideal in $T$.  Suppose that $I$ can be generated by a $d$-sequence. Then: 

{\rm (a)[Herzog-Simis-Vasconcelos, K\"uhl, Ng\^o]} We have $\reg(\Rees(I))=0$.

{\rm (b)[Huneke]} The natural surjection of $T$-algebras 
$\oplus_{k\geq 0}\Sym^k(I)\to\Rees(\Fp)$ is an isomorphism. 

{\rm (c)[Huneke]} If $I$ can be generated by a Cohen-Macaulay $d$-sequence then the natural surjection of $T$-algebras 
\mbox{$\oplus_{k\geq 0}\Sym^k(I/I^2)\to\oplus_{k\geq 0}I^k/I^{k+1}$} is an isomorphism.

{\rm (d)[Huneke]} In the situation of Proposition \ref{propH0}, the ideal 
$\Fp$ is generated by a Cohen-Macaulay $d$-sequence iff $T/\Fp$ is a Cohen-Macaulay ring.
\label{theorVR}
\end{theor}
\beginProof 
See \cite[Cor. 1.2]{Trung-C} for (a) (and for more references - this is a deep result). 
See \cite[Th. 3.1]{Huneke-Sym}  for (b). See \cite[Theorem C]{Huneke-Koszul} for (c) and 
\cite[Prop. 2.4]{Huneke-Koszul} for (d). 
\endProof

We now have

\begin{prop}Assumptions as in Corollary \ref{corFPG}. Suppose also that $X$ is integral and that $X_G$ is a 
disjoint union of integral schemes. Suppose that $X$ is regular at the generic points of the irreducible components of $X_G.$ 
Suppose furthermore that for each $x\in X_G$, the local ring $\CO_{X,x}$ of $x$ in $X$ is Cohen-Macaulay with infinite residue field. Finally, suppose the for any $x\in X_G$, the $\CO_{X,x}$-module $(N_{X_G/X})_{\CO_{X,x}}$ can be generated by 
${\rm height}(I\cdot\CO_{X,x})+1$ (or less) elements. 

Then the conclusion of Corollary \ref{corFPG} holds for $\lambda=0$. 

Furthermore, if $X_G$ is also Cohen-Macaulay, we have 
$I^u/I^{u+1}\simeq\Sym^{u}(I/I^2)$ for all $u\geq 0$. 
\label{propREDL}
\end{prop}
Note the following. Let $C$ be an irreducible component of $X_G$ and let $\eta_C$ be its generic point. Then  ${\rm height}(I\cdot\CO_{X,x})$ is constant on $C$, with value ${\rm dim}(\CO_{X,\eta_C})={\rm codim}(C,X)$.

\beginProof Nakayama's lemma implies that 
$I\cdot \CO_{X,x}$ can be generated by ${\rm height}(I\cdot\CO_{X,x})+1$  elements. Hence the conclusion of the proposition follows from Proposition \ref{propH0} and Theorem \ref{theorVR} (a), (c), (d).\endProof

\begin{ex}\rm Here is a simple example. Suppose that $S$ is the spectrum of 
an infinite perfect field and that $Y$ is a geometrically integral plane projective curve. Then $\Omega_{Y/S}$ is generated 
locally by two elements, since we have an epimorphism of sheaves $\Omega_{\mP^2/S}|Y\to \Omega_{Y/S}.$ Now set $X:=Y\times_S Y$ and let $G$ act on $X$ by swapping the factors. 
Then $Y\simeq X_G$ and so $X_G$ is integral and has codimension one in $X$. Furthermore, 
$N_{X_G/X}\simeq\Omega_{Y/S}$ can be locally generated by $\cod(X_G,X)+1=2$ elements. The local rings of $Y$ are Cohen-Macaulay with infinite residue fields since $Y$ is lci over $S$ and so the local rings of $X$ are also Cohen-Macaulay with infinite residue fields. 
Furthermore, $X$ is regular at the generic point of $Y$, since $Y$ is geometrically integral over $S$. 
Hence $X$ satisfies the assumptions of Proposition \ref{propREDL}. 
This example will be generalised in Proposition \ref{propSLCI} below, to whose proof we refer for more details.\end{ex}

\section{The Riemann-Roch theorem}
\label{secMT}

We use the notation of the introduction. We recall some of it for the convenience of the reader.

Let $S$ be a separated noetherian scheme on which $2$ is invertible and let $f:Y\to S$ be a proper and flat scheme over 
$S$. Let $V$ be a vector bundle on $Y$ and suppose that the coherent sheaves 
$\R^i f_*(V)$ are all locally free. Let $X:=Y\times_S Y$ and let $\pi:\wt{X}\to Y\times_S Y$ be 
the blow-up of $Y\times_S Y$ along the diagonal $\Delta$ of $Y\times_S Y$. Let $\phi:E\to\Delta$ be the corresponding 
exceptional divisor. Let $I_\Delta$ be the ideal sheaf of $\Delta$.  

Let $\cm(f)\geq 0$ be the minimal natural number $\lambda$ such that 
$\R^a\pi_*(\CO(-E)^{\otimes r})=0$ for all $a>0$ and such that the natural 
morphism of sheaves $I_\Delta^r\to \pi_*(\CO(-E)^{\otimes r})$ is an isomorphism  for all $r\geq \lambda.$  As already mentioned, this exists by \cite[Cohomology of Schemes, Lemmas 14.2 and 14.3]{StacksProject}. Note that we then have $\R^\bullet\phi(\CO(-E)|_E^{\otimes r})=(I^r_\Delta/I^{r+1}_\Delta)^r$ 
(in the derived category) for all $r\geq \lambda.$ To see this, apply $\R^\bullet\pi_*(\cdot)$
to the exact sequence
$$
0\to\CO(-(r+1)E)\to \CO(-rE)\to \CO(-E)|_E^{\otimes r}\to 0.
$$
Note that if 
$Y$ is smooth over $S$ (so that the diagonal immersion is a regular immersion) then  $\cm(f)=0$. 

Let 
$\delta_0(f):=\delta_0(N_{E/\wt{X}})$ be the smallest natural number $n$ such that 
\mbox{$(N_{E/\wt{X}}-1)^{\otimes(n +1)}=0$} in $\r K^0(E)_\mQ$. Here $\r K^0(E)_\mQ$ is 
the quotient of $K^0(E)_\mQ$ by the annihilator of $K_0(E)_\mQ.$

Finally, for any $\lambda,\delta\geq 0$ let 
$$
\GTI(f,\delta,\lambda):={(-1)^\lambda\over 2}\sum_{j=0}^\delta E_j \sum_{k=0}^{\delta}\sum_{u=0}^k\ {(-1)^{k-u}s(k,j)\over u!(k-u)!}I^{u+\lambda}_\Delta/I^{u+\lambda+1}_\Delta+\sum_{k=0}^{\lambda-1}(-1)^k I^k_\Delta/I^{k+1}_\Delta\in K_0(Y)_\mQ.
$$ 

We now consider $X=Y\times_S Y$ as a $G$-equivariant scheme, with the action of $G=(\mZ/2\mZ)_S$ 
which swaps the two factors. We then have a natural identification $X_G=\Delta\simeq Y$. 

\begin{theor} 
The local term of $X$ is given by the image of $\GTI(f,\delta,\lambda)$ in \mbox{$K_0(X_G,G)\otimes_{R(G)}\mQ$} for 
any $\delta\geq\delta_0(f)$ and any $\lambda\geq \cm(f)$.
\label{theorLTP}
\end{theor}
\beginProof Note that 
any point of $\Delta$  has an open affine $G$-invariant neighbourhood, which lies over an open affine subscheme of $S.$ To see this, let $x\in\Delta$ and view $x$ as a point of $Y$. 
Then this point has an open affine neighbourhood $U$ which lies over an open affine subscheme of $S$. Hence $x\in U\times_S U$ in $X$ and $U\times_S U$ is affine, $G$-invariant and lies 
over an open affine subscheme of $S$. So the assumption of Corollary \ref{corFPG} are satisfied and the conclusion follows from the conclusion of Corollary \ref{corFPG}.\endProof

In particular, $\GTI(f,\delta,\lambda)$ is constant in the range $\delta\geq\delta_0(f), \lambda\geq \cm(f).$ We write $\GTI(f)$ for this constant value.
 
\begin{cor}[Theorem \ref{theorMT}] 
Suppose that $V$ is a vector bundle on $Y$ and 
that $\R^i f_*(V)$ is locally free for all $i\geq 0$. Then 
the equality 
\begin{eqnarray}
&&\sum_i(-1)^i[\Sym^2(\R^i f_*(V)-\Lambda^2(\R^i f_*(V))]=\R^\bullet f_*(\GTI(f)\otimes(\Sym^2(V)-\Lambda^2(V)))
\end{eqnarray}
holds in $K_0(S)_\mQ.$
\label{theorMTP}
\end{cor}

\beginProof (of Corollary \ref{theorMTP}). 
We let $W:=\pi_1^*(V)\otimes\pi_2^*(V)$, where $\pi_{1,2}:Y\times_S Y\to Y$ are the natural projections.
The bundle $W$ is naturally endowed with the $G$-action which swaps the factors of the tensor product. 
Hence we deduce from Theorem \ref{theorLTP} and the projection formula that
\begin{eqnarray*}
&&\Tr(\R^\bullet (f\times f)_*(W))=\Tr(\R^\bullet f_*(\GTI(f,\delta,\lambda)\otimes W|_Y)\\&=&
\Tr(\R^\bullet f_*(\Big[{(-1)^\lambda\over 2}\sum_{j=0}^\delta E_j \sum_{k=0}^{\delta}\sum_{u=0}^k\ {(-1)^{k-u}s(k,j)\over u!(k-u)!}[I^{u+\lambda}_\Delta/I^{u+\lambda+1}_\Delta]_\triv\\&+&\sum_{k=0}^{\lambda-1}(-1)^k [I^k_\Delta/I^{k+1}_\Delta]_\triv\Big]\otimes W|_Y))
\end{eqnarray*}
in $K_0(S)_\mQ$ for any $\delta\geq\delta_0(f)$ and any $\lambda\geq \cm(f)$. 
Now note that by definition $(W|_G)_{+}=\Sym^2(V)$ and $(W|_G)_{-}=\Lambda^2(V).$ Hence we get the equality
\begin{eqnarray}
&&\Tr(\R^\bullet (f\times f)_*(W))=\R^\bullet f_*(\GTI(f,\delta,\lambda)\otimes(\Sym^2(V)-\Lambda(V))\nonumber\\&=&
\R^\bullet f_*(\Big[{(-1)^\lambda\over 2}\sum_{j=0}^\delta E_j \sum_{k=0}^{\delta}\sum_{u=0}^k\ {(-1)^{k-u}s(k,j)\over u!(k-u)!}I^{u+\lambda}_\Delta/I^{u+\lambda+1}_\Delta\nonumber\\&+&\sum_{k=0}^{\lambda-1}(-1)^k I^k_\Delta/I^{k+1}_\Delta\Big]\otimes (\Sym^2(V)-\Lambda(V)))\label{eqYAFE}
\end{eqnarray}
in $K_0(S)_\mQ$. Here the two expressions on the right are computed in $K_0(X)_\mQ$ without 
reference to any equivariant structure.

The corollary now follows from \refeq{eqYAFE} and the next 

\begin{sublem}We have 
$$
\Tr(\R^\bullet (f\times f)_*(W))=\sum_i(-1)^i[\Sym^2(\R^i f_*(V)-\Lambda^2(\R^i f_*(V))]
$$
in $K_0(S)_\mQ$.  
\label{slemCRUX}
\end{sublem}

\beginProof (of sublemma \ref{slemCRUX}). 
By the K\"unneth formula (see \cite[III, par. 6, Th. 6.7.3]{EGA}), we have a canonical isomorphism 
\begin{equation}
R^i (f\times f)_*(W)=\R^i (f\times f)_*(\pi_1^*(V)\otimes\pi_2^*(V))\cong\bigoplus_{t}\R^t f_*(V)\otimes\R^{i-t}f_*(V).
\label{kis}
\end{equation}
The vector bundle $$\bigoplus_{t}\R^t f_*(V)\otimes\R^{i-t}f_*(V)$$ carries a natural $G$-action by permutation, namely the 
action such that the $1\in\mZ/2\mZ$ sends
$
\bigoplus_t w_t\otimes w_{i-t}
$
to
\mbox{$
\bigoplus_t (-1)^{t(i-t)}w_{i-t}\otimes w_{t}.
$}
 By the Koszul rule of signs, the isomorphism \refeq{kis} becomes $G$-equivariant with this choice of $G$-action on the righthand side.  

Now let us first suppose that $i$ is odd. We then have an equivariant isomorphism
\begin{eqnarray}
\R^i (f\times f)_*(\pi_1^*(V)\otimes\pi_2^*(V))\cong\bigoplus_{0\leq t\leq \lfloor i/2\rfloor}\R^t f_*(V)\otimes\R^{i-t}f_*(V)\oplus 
\R^{i-t} f_*(V)\otimes\R^{t}f_*(V).\label{kiseq}
\end{eqnarray}
where  (by the above) the element  $1\in\mZ/2\mZ$ acts on  
$\R^t f_*(V)\otimes\R^{i-t}f_*(V)\oplus 
\R^{i-t} f_*(V)\otimes\R^{t}f_*(V)$  by the automorphism which swaps the two factors of the direct sum. We thus have non equivariant isomorphisms
\begin{eqnarray*}
&&\big(\R^t f_*(V)\otimes\R^{i-t}f_*(V)\oplus 
\R^{i-t} f_*(V)\otimes\R^{t}f_*(V)\big)_{+}\\&\simeq&\R^t f_*(V)\otimes\R^{i-t}f_*(V)\simeq 
\big(\R^t f_*(V)\otimes\R^{i-t}f_*(V)\oplus 
\R^{i-t} f_*(V)\otimes\R^{t}f_*(V)\big)_{-}
\end{eqnarray*}
and hence $\Tr(\R^i (f\times f)_*(\pi_1^*(V)\otimes\pi_2^*(V))=0$. 

Now suppose that $i$ is even. We then have an equivariant isomorphism
\begin{eqnarray*}
&&\R^i (f\times f)_*(\pi_1^*(V)\otimes\pi_2^*(V))\\&\cong&\R^{i/2} f_*(V)\otimes\R^{i/2}f_*(V)\oplus\bigoplus_{0\leq t<i/2}\big(\R^t f_*(V)\otimes\R^{i-t}f_*(V)\oplus 
\R^{i-t} f_*(V)\otimes\R^{t}f_*(V)\big).\label{kiseqp}
\end{eqnarray*}
Here the element  $1\in\mZ/2\mZ$ acts on  $\R^t f_*(V)\otimes\R^{i-t}f_*(V)\oplus 
\R^{i-t} f_*(V)\otimes\R^{t}f_*(V)$ by multiplication by $(-1)^t$ followed by the automorphism which swaps the two factors of the direct sum. On the other hand, the element  $1\in\mZ/2\mZ$ acts on $\R^{i/2} f_*(V)\otimes\R^{i/2}f_*(V)$ by multiplication by $(-1)^{i/2}.$

So we conclude that
\begin{eqnarray*}
&&\Tr(\R^i (f\times f)_*(\pi_1^*(V)\otimes\pi_2^*(V)))\\&=&(-1)^{i/2}\big((\R^{i/2} f_*(V)\otimes\R^{i/2}f_*(V))_{+}-(\R^{i/2} f_*(V)\otimes\R^{i/2}f_*(V))_{-}\big)\\&=&(-1)^{i/2}\big(\Sym^2(\R^{i/2} f_*(V))-
\Lambda^2(\R^{i/2} f_*(V))\big)
\end{eqnarray*}
Summing over all $i\geq 0$ we get the equality of the sublemma. 
\endProof

\section{Special cases}

\label{secSC}

\subsection{local complete intersections}

We shall now consider Corollary \ref{theorMTP} in the situation where the morphism $f$ is an lci morphism. The following Proposition in particular proves the ARR theorem for any 
flat, lci and projective scheme over a noetherian separated scheme, which carries an 
ample line bundle and over which $2$ is invertible.

\begin{prop}
Assumption as in Corollary \ref{theorMTP}. Suppose  that there exists 
a regular closed $S$-immersion $Y\hookrightarrow W$, where $W$ is smooth over $S$. Suppose also 
that $Y$ carries an ample family of line bundles. Then we have
\begin{eqnarray*}
\GTI(f)
&=&\Theta^2(N_{Y/W})\otimes\Theta^2(\Omega_{W/S}|_Y)^{-1}
\end{eqnarray*}
in $K_0(Y)_\mQ$. 
\label{propLCIT}
\end{prop}
\beginProof We again consider $Y\times_S Y$ and $W\times_S W$ as $G$-equivariant schemes. The natural embedding $Y\times_S Y\hookrightarrow W\times_S W$ is then $G$-equivariant 
and the scheme $W\times_S W$ is smooth over $S$.  As before, we have 
$W\simeq (W\times_S W)_G$ and $Y\simeq (Y\times_S Y)_G$ via the diagonal immersions. 

From Theorem \ref{theorLTP} and Corollary \ref{corGOODLT} we get 
\begin{equation}
\GTI(f)=\big(\sum_{i\geq 0}(-1)^i\Tr(\rTor^i_{\CO_{W\times_S W}}(\CO_{W},\CO_{Y\times_S Y}))\big)\otimes_{\CO_Y}\Theta^2(N_{W/W\times_S W}|_{Y})^{-1}
\label{eqGTILCI}
\end{equation}
in $K_0(Y)_\mQ$ (here $N_{W/W\times_S W}$ is not considered as an equivariant sheaf, but $\rTor^i_{\CO_Z}(\CO_{Z_G},\CO_X)$ carries its natural equivariant structure). Note that $N_{W/W\times_S W}\simeq\Omega_{W/S}$ by definition. 

We shall now compute 
$\rTor^i_{\CO_{W\times_S W}}(\CO_{W},\CO_{Y\times_S Y}).$ Let 
$$
V^\bullet:\dots\to V_l\to V_{l-1}\to \dots V_0\to\CO_Y\to 0
$$
be a (possibly infinite) resolution of $\CO_Y$ in $W$ by flat quasi-coherent sheaves. 
This exists by \cite[Appendix A.1]{Efim-Coh}.  Let $\pi_1,\pi_2:Z\to W$ 
be the two natural projections. Since 
$W$ is flat over $S$, the sequence $\pi_1^*(V^\bullet)\otimes\pi_1^*(V^\bullet)$ is 
a resolution of $Y\times_S Y$ in $W\times_S W$. By construction, this resolution carries 
a $G$-equivariant structure (the $G$-action permuting the factors according to the Koszul rule of signs) making it into an equivariant resolution of $Y\times_S Y$. By definition, we have an equivariant 
isomorphism of sheaves
$$
\CH^{-i}((\pi_1^*(V^\bullet)\otimes\pi_1^*(V^\bullet))|_{W})\simeq \rTor^i_{\CO_{W\times_S W}}(\CO_{W},\CO_{Y\times_S Y})
$$
Furthermore, from the definitions we see that we have 
$$
\pi_1^*(V^\bullet)|_{W}=\pi_2^*(V^\bullet)|_{W}=V^\bullet
$$
and thus 
$$
\CH^{-i}((\pi_1^*(V^\bullet)\otimes \pi_1^*(V^\bullet))|_{W})\simeq
\CH^{-i}(V^\bullet\otimes V^\bullet)\simeq 
\rTor^i_{\CO_{W}}(\CO_{Y},\CO_Y).
$$
Since $Y$ is regularly immersed in $W$, we have a functorial isomorphism 
$$
\rTor^i_{\CO_{W}}(\CO_{Y},\CO_Y)\simeq\Lambda^i(N_{Y/W})
$$
(see \cite[VII, Prop 2.5 i')]{SGA6}). Finally, since $1\in\mZ/2\mZ$ acts by $(-1)^i$ on $\Lambda^i(N_{Y/W})$, we 
get from \refeq{eqGTILCI} that
$$
\GTI(f)=\Theta^2(N_{Y/W})\otimes\Theta^2(\Omega_{W/S}|_Y)^{-1}
$$
as required.\endProof

\begin{rem}\rm Notice that if $S$ is regular then $W$ carries an ample family of line bundles and 
then so does $Y$ (see \cite[II 2.2.7.1]{SGA6}). In that case, we have $K_0(S)\simeq K^0(S)$ and 
Proposition \ref{propLCIT} is stronger than the classical GRR theorem, because 
no projectivity assumptions are made.
\end{rem}

\subsection{Cohen-Macaulay schemes with hypersurface singularities}

If $W$ is a noetherian scheme, we shall say that $W$ has at most hypersurface singularities 
if for any $w\in W$, the embedding dimension of $\CO_{W,w}$ is at most one more than 
the dimension of $\CO_{W,w}$. 

For the proof of the next proposition, we shall need the

\begin{lemma}
Let $K$ be field and let $Z$ be a scheme of finite type over $K$. 
Suppose that the geometric fibre of $Z$ over $K$ is integral with  at most hypersurface singularities. 
Then for every point $z\in Z$, the $\CO_{Z,z}$-module $\Omega_{Z,\CO_{Z,z}}$ is generated 
by  $\dim(Z)+1$ (or less) elements.
\label{lemINT}
\end{lemma}
\beginProof Let $n\geq 0$ and let 
$$
Z_n:=\{z\in Z\,|\,\rk_{\kappa(z)}(\Omega_{Z,\kappa(z)})\geq n\}
$$
where $\kappa(z)$ is the residue field at $z$. This set is closed by \cite[II.5, Exercise 5.8 a)]{Hartshorne-AG} and we view 
it as a reduced closed subscheme. By Nakayama's lemma, the conclusion of the lemma is equivalent to saying that 
$Z_n=\emptyset$ if \mbox{$n>\dim(Z)+1$.} So suppose for contradiction that we are given $n>\dim(Z)+1$ 
and that $Z_n\not=\emptyset.$ Suppose that $n$ is maximal with that property. Let $\iota:Z_n\to Z$ be the inclusion morphism. By \cite[II.5, Exercise 5.8 c)]{Hartshorne-AG}, the sheaf $\iota^*(\Omega_{Z/K})$ is then locally free of rank $n$. 
In particular, for any closed point $\bar z$ of $Y_{\bar K}$ lying over a  point of $Z_n$, a  
minimal set of generators of $\Omega_{Z_{\bar K},\CO_{Z_{\bar K},\bar z}}$ has cardinality 
$n$. On the other hand, $\dim(Z_{\bar K})=\dim(Z)$ and $\dim(\CO_{Z_{\bar K},\bar z})=\dim(Z_{\bar K})$ 
since $\bar z$ is closed. Finally, we have $\Fm_{\bar z}/\Fm_{\bar z}^2\simeq\Omega_{Z_{\bar K},\kappa(\bar z)}$, again because $\bar z$ is closed. So
 $n=\rk_{\kappa(\bar z)}(\Fm_{\bar z}/\Fm_{\bar z}^2)>\dim(\CO_{Z_{\bar K},\bar z})+1$, which contradicts 
 the fact that $Z_{\bar K}$ has only hypersurface singularities.\endProof

\begin{prop} Assumptions as in Corollary \ref{theorMTP}. 
Suppose that $S$ is Cohen-Macaulay, integral with infinite residue fields, and that the regular locus $\reg(S)$ of $S$ contains 
a non-empty open set. Suppose that the geometric fibres of $Y$ over $S$ are integral, Cohen-Macaulay and 
have at most hypersurface singularities.   Then for all 
$\delta\geq\delta_0(f)$ we have
\begin{eqnarray*}
&&\GTI(f)={1\over 2}\sum_{j=0}^\delta E_j\sum_{k=0}^{\delta}\sum_{u=0}^k\ {(-1)^{k-u}s(k,j)\over (k-u)!u!}\Sym^u(\Omega_{Y/S}).
\end{eqnarray*}
\label{propSLCI}
 Furthermore, we have $\Sym^u(\Omega_{Y/S})\simeq I_\Delta^u/I_\Delta^{u+1}$ for all 
$u\geq 0.$
\end{prop}
Note that $\Omega_{Y/S}$ might not be locally free. We also note that the assumptions 
imply that $Y$ is integral. See the proof below.
\beginProof 
By \cite[Prop. 4.3.8]{Liu-Alg}, the schemes  $Y$ and $Y\times_S Y$ are integral. By \cite[tag045J, Lemma 10.163.3]{StacksProject}, $Y$ is also Cohen-Macaulay and so is $Y\times_S Y$. Furthermore, 
from the assumptions we see that $Y\times_S Y$ is regular at the generic point of $Y$. 

Let now $y\in Y$. 
By  Lemma \ref{lemINT} and Nakayama's lemma 
the $\CO_{Y_{f(y),y}}$-module $\Omega_{Y/S,\CO_{Y_{f(y),y}}}$ is generated by $d_y+1$ elements, 
where $d_y=\dim(Y_{f(y)}).$ Now note that $d_y$ does not depend on $y$ by \cite[tag02NI, Lemma 29.29.4]{StacksProject} (because $f$ is flat with Cohen-Macaulay fibres).  Let 
$\eta$ (resp. $\rho$) be the generic point of $Y$ (resp. $S$).  Since $f$ is flat and $S$ is integral, the point 
$\eta$ lies over $\rho$. Therefore, by \cite[4.3.12]{Liu-Alg}, we have 
$\dim(\CO_{Y\times_S Y,\eta})=\dim(\CO_{Y_\rho\times_\rho Y_\rho,\eta})$ and since 
$Y_\rho\times_\rho Y_\rho$ is an  integral scheme of finite type over a field, we have 
$$\dim(\CO_{Y_\rho\times_\rho Y_\rho,\eta})={\rm codim}(Y_\rho,Y_\rho\times_\rho Y_\rho)=
\dim(Y_\rho)=d_y.$$ We conclude that the $\CO_{Y_{f(y),y}}$-module $\Omega_{Y/S,\CO_{Y_{f(y),y}}}$is generated 
by $\dim(\CO_{Y\times_S Y,\eta})+1$ elements. Hence by Proposition \ref{propREDL} we have $\cm(f)=0$ 
 and $\Sym^u(\Omega_{Y/S})\simeq I_\Delta^u/I_\Delta^{u+1}$. 
The formula of the proposition now follows from Corollary \ref{theorMTP}.\endProof 

\begin{rem}\rm  We keep the assumptions of Proposition \ref{propSLCI}. Suppose in addition that $Y$ is a Cartier divisor in an $S$-scheme $H$, which is smooth over $S$, and that $Y$ carries 
an ample family of line bundles. Then Proposition \ref{propSLCI} applies but 
so does Proposition \ref{propLCIT}. So suppose for the time of this remark that we are in that situation. We then have a resolution
$$
0\to N_{Y/H}\to \Omega_{H/S}|_Y\to \Omega_{Y/S}\to 0
$$
 In particular $\Omega_{Y/S}$ represents a strictly   perfect complex in the derived category of $\CO_Y$-modules, namely the cotangent complex of $Y$ over $S$. Also, since $N_{Y/S}$ is a line bundle and 
 $Y$ is smooth over $S$ at its generic point, we have an exact sequence
 $$
 0\to N_{Y/H}\otimes\Sym^{n-1}(\Omega_{H/S}|_Y)\to\Sym^{n}(\Omega_{H/S}|_Y)\to\Sym^n(\Omega_{Y/S})\to 0
 $$
 for all $n\geq 0$ (see \cite[tag01CF, Lemma 17.21.4]{StacksProject}) so that the coherent sheaf $\Sym^u(\Omega_{Y/S})$ represents the $u$-th symmetric power of 
 the cotangent complex of $Y$ over $S$. We will give a direct proof (not involving the fixed point formula) of the fact that 
 the right-hand side of Propositions \ref{propSLCI} and \ref{propLCIT} coincide (note that we now a priori that they coincide since they are both equal to $\GTI(f)$). This will show that under the given assumptions 
Proposition \ref{propSLCI} is a natural generalisation of the ARR theorem to a non-projective situation. 

We need to show that 
\begin{eqnarray}
&&\Theta^2(\Omega_{H/S}|_Y)^{-1}\otimes \Theta^2(N_{Y/H})\nonumber\\&=&
{1\over 2}\sum_{j=0}^\delta E_j \sum_{k=0}^{\delta}\sum_{u=0}^k\ {(-1)^{k-u}s(k,j)\over (k-u)!u!}
\big(\Sym^{u}(\Omega_{H/S}|_Y)-N_{Y/H}\otimes\Sym^{u-1}(\Omega_{H/S}|_Y)\big)\label{eqBBFIN}
\end{eqnarray}
in $K_0(Y)_\mQ$ for any $\delta\geq\delta_0(f).$ Here the inverse of $\Theta^2(\Omega_{H/S}|_Y)$
 is taken in $K^0(Y)_\mQ$ and we take the image of this inverse in $K_0(Y)_\mQ$.

Let $P(x)=\sum_j a_j x^j\in K^0(X)_\mQ[x]$ be a polynomial such that 
$P(k)=\Sym^k(\Omega_{H/S}|_Y)$ in $K^0(Y)_\mQ$. This exists by Proposition \ref{propKEYI} (b) (1). 
We then have
$$
\Sym^n(\Omega_{H/S}|_Y)-N_{Y/H}\otimes\Sym^{n-1}(\Omega_{H/S}|_Y)=P(n)-N_{Y/H}\otimes P(n-1)
$$
(note that one can show that $P(-1)=0$ [we skip the proof] so this makes sense) and we let 
$$
Q(x):=P(x)-N_{Y/H}\otimes P(x-1)=:\sum_j b_j x^j\in K^0(Y)_\mQ[x].
$$ We  have by construction an identity of power 
series 
$$
\Lambda_{-t}(N_{Y/H})\Sym_t(\Omega_{H/S}|_Y)=\sum_k Q(n) t^n\in K^0(Y)_\mQ[[t]]
$$
and applying Proposition \ref{propKEYI} (a) \& (b) (2) we get 
\begin{equation}
\Theta^2(\Omega_{H/S}|_Y)^{-1}\otimes \Theta^2(N_{Y/H})=
{1\over 2}\sum_j E_j b_j\label{eqBFIN}
\end{equation}
in $K^0(Y)_\mQ$. On the other hand, for all $n\geq 0$ we have 
$$
\R^\bullet \phi_*(\CO(n))=\Sym^n(\Omega_{Y/S})=\Sym^n(\Omega_{H/S}|_Y)-N_{Y/H}\otimes\Sym^{n-1}(\Omega_{H/S}|_Y)
$$
in $K_0(Y)_\mQ$ and so by Lemma \ref{lemPOL} (a) \& (c) (1), we have 
$$
\Sym^n(\Omega_{Y/S})=\sum_{j=0}^{\delta}\big[\sum_{k=0}^{\delta}\sum_{u=0}^k\ {(-1)^{k-u}s(k,j)\over u!(k-u)!}\Sym^u(\Omega_{Y/S})\big] n^j
$$
in $K_0(Y)_\mQ$ for $\delta\geq\delta_0(f)$. In particular, we have 
$$Q(n)=\sum_{j=0}^{\delta}\big[\sum_{k=0}^{\delta}\sum_{u=0}^k\ {(-1)^{k-u}s(k,j)\over u!(k-u)!}\Sym^u(\Omega_{Y/S})\big] n^j$$
in $K_0(Y)_\mQ$ (not $K^0(Y)_\mQ$ !) for all $n\geq 0$. By Lemma \ref{lemVS}, we see that the image of $b_j$ in $K_0(Y)_\mQ$ 
is $$\sum_{k=0}^{\delta}\sum_{u=0}^k\ {(-1)^{k-u}s(k,j)\over u!(k-u)!}\Sym^u(\Omega_{Y/S})=
\sum_{k=0}^{\delta}\sum_{u=0}^k\ {(-1)^{k-u}s(k,j)\over u!(k-u)!}
\big(\Sym^u(\Omega_{H/S}|_Y)-N_{Y/H}\otimes\Sym^{u-1}\big).$$
Combining this with \refeq{eqBFIN} we get \refeq{eqBBFIN}. 
\end{rem}


\begin{bibdiv}
\begin{biblist}

\bib{BFM-RR}{article}{
   author={Baum, Paul},
   author={Fulton, William},
   author={MacPherson, Robert},
   title={Riemann-Roch for singular varieties},
   journal={Inst. Hautes \'{E}tudes Sci. Publ. Math.},
   number={45},
   date={1975},
   pages={101--145},
   issn={0073-8301},
}

\bib{BGG}{article}{
   author={Be\u{\i}linson, A. A.},
   title={The derived category of coherent sheaves on ${\bf P}^n$},
   note={Selected translations},
   journal={Selecta Math. Soviet.},
   volume={3},
   date={1983/84},
   number={3},
   pages={233--237},
   issn={0272-9903},
}

\bib{Brod-LC}{book}{
   author={Brodmann, M. P.},
   author={Sharp, R. Y.},
   title={Local cohomology},
   series={Cambridge Studies in Advanced Mathematics},
   volume={136},
   edition={2},
   note={An algebraic introduction with geometric applications},
   publisher={Cambridge University Press, Cambridge},
   date={2013},
   pages={xxii+491},
   isbn={978-0-521-51363-0},
}

\bib{Efim-Coh}{article}{
   author={Efimov, Alexander I.},
   author={Positselski, Leonid},
   title={Coherent analogues of matrix factorizations and relative
   singularity categories},
   journal={Algebra Number Theory},
   volume={9},
   date={2015},
   number={5},
   pages={1159--1292},
   issn={1937-0652},
   doi={10.2140/ant.2015.9.1159},
}

\bib{Eis-Co}{book}{
   author={Eisenbud, David},
   title={Commutative algebra},
   series={Graduate Texts in Mathematics},
   volume={150},
   note={With a view toward algebraic geometry},
   publisher={Springer-Verlag, New York},
   date={1995},
   pages={xvi+785},
   isbn={0-387-94268-8},
   isbn={0-387-94269-6},
   doi={10.1007/978-1-4612-5350-1},
}




\bib{Fresan-Periods}{article}{
   author={Fres\'{a}n, Javier},
   title={Periods of Hodge structures and special values of the gamma
   function},
   journal={Invent. Math.},
   volume={208},
   date={2017},
   number={1},
   pages={247--282},
   issn={0020-9910},
   doi={10.1007/s00222-016-0690-4},
}

\bib{FL-RR}{book}{
   author={Fulton, William},
   author={Lang, Serge},
   title={Riemann-Roch algebra},
   series={Grundlehren der mathematischen Wissenschaften [Fundamental
   Principles of Mathematical Sciences]},
   volume={277},
   publisher={Springer-Verlag, New York},
   date={1985},
   pages={x+203},
   isbn={0-387-96086-4},
   doi={10.1007/978-1-4757-1858-4},
}

\bib{SGA6}{collection}{
   title={Th\'{e}orie des intersections et th\'{e}or\`eme de Riemann-Roch},
   language={French},
   series={Lecture Notes in Mathematics, Vol. 225},
   note={S\'{e}minaire de G\'{e}om\'{e}trie Alg\'{e}brique du Bois-Marie 1966--1967 (SGA 6);
   Dirig\'{e} par P. Berthelot, A. Grothendieck et L. Illusie. Avec la
   collaboration de D. Ferrand, J. P. Jouanolou, O. Jussila, S. Kleiman, M.
   Raynaud et J. P. Serre},
   publisher={Springer-Verlag, Berlin-New York},
   date={1971},
   pages={xii+700},
}

\bib{EGA}{article}{
   author={Grothendieck, A.},
   author={Dieudonn\'e, J.},
   title={{\'El\'ements de g\'eom\'etrie alg\'ebrique.} 
 {\rm Inst. Hautes \'Etudes Sci. Publ. Math.} {\bf 4, 8, 11, 17, 20, 
24, 28, 32} {\rm (1960-1967)}.}
}

\bib{Hartshorne-AG}{book}{
   author={Hartshorne, Robin},
   title={Algebraic geometry},
   series={Graduate Texts in Mathematics, No. 52},
   publisher={Springer-Verlag, New York-Heidelberg},
   date={1977},
   pages={xvi+496},
   isbn={0-387-90244-9},
}

\bib{Hartshorne-Residues}{book}{
   author={Hartshorne, Robin},
   title={Residues and duality},
   series={Lecture Notes in Mathematics, No. 20},
   note={Lecture notes of a seminar on the work of A. Grothendieck, given at
   Harvard 1963/64;
   With an appendix by P. Deligne},
   publisher={Springer-Verlag, Berlin-New York},
   date={1966},
   pages={vii+423},
}

\bib{HSV-ACI}{article}{
   author={Hong, Jooyoun},
   author={Simis, Aron},
   author={Vasconcelos, Wolmer V.},
   title={The equations of almost complete intersections},
   journal={Bull. Braz. Math. Soc. (N.S.)},
   volume={43},
   date={2012},
   number={2},
   pages={171--199},
   issn={1678-7544},
   doi={10.1007/s00574-012-0009-z},
}

\bib{Huneke-d}{article}{
   author={Huneke, Craig},
   title={The theory of $d$-sequences and powers of ideals},
   journal={Adv. in Math.},
   volume={46},
   date={1982},
   number={3},
   pages={249--279},
   issn={0001-8708},
   doi={10.1016/0001-8708(82)90045-7},
}


\bib{Huneke-Koszul}{article}{
   author={Huneke, Craig},
   title={The Koszul homology of an ideal},
   journal={Adv. in Math.},
   volume={56},
   date={1985},
   number={3},
   pages={295--318},
   issn={0001-8708},
   doi={10.1016/0001-8708(85)90037-4},
}

\bib{Huneke-Sym}{article}{
   author={Huneke, Craig},
   title={On the symmetric and Rees algebra of an ideal generated by a
   $d$-sequence},
   journal={J. Algebra},
   volume={62},
   date={1980},
   number={2},
   pages={268--275},
   issn={0021-8693},
   doi={10.1016/0021-8693(80)90179-9},
}

\bib{Koeck-KC}{article}{
   author={K\"{o}ck, Bernhard},
   title={Computing the homology of Koszul complexes},
   journal={Trans. Amer. Math. Soc.},
   volume={353},
   date={2001},
   number={8},
   pages={3115--3147},
   issn={0002-9947},
   doi={10.1090/S0002-9947-01-02723-4},
}

\bib{Kim-Alt}{article}{
   author={Kim, Taekyun},
   author={Rim, Seog-Hoon},
   author={Simsek, Yilmaz},
   title={A note on the alternating sums of powers of consecutive
   $q$-integers},
   journal={Adv. Stud. Contemp. Math. (Kyungshang)},
   volume={13},
   date={2006},
   number={2},
   pages={159--164},
   issn={1229-3067},
}

\bib{Liu-Alg}{book}{
   author={Liu, Qing},
   title={Algebraic geometry and arithmetic curves},
   series={Oxford Graduate Texts in Mathematics},
   volume={6},
   note={Translated from the French by Reinie Ern\'{e};
   Oxford Science Publications},
   publisher={Oxford University Press, Oxford},
   date={2002},
   pages={xvi+576},
   isbn={0-19-850284-2},
}



\bib{Mat-Co}{book}{
   author={Matsumura, Hideyuki},
   title={Commutative ring theory},
   series={Cambridge Studies in Advanced Mathematics},
   volume={8},
   edition={2},
   note={Translated from the Japanese by M. Reid},
   publisher={Cambridge University Press, Cambridge},
   date={1989},
   pages={xiv+320},
   isbn={0-521-36764-6},
}

\bib{Magma}{misc}{
title={Magma calculator. Many different contributors.},
note={\url{http://magma.maths.usyd.edu.au/calc/}}
}


\bib{Nori-HRR}{article}{
   author={Nori, Madhav V.},
   title={The Hirzebruch-Riemann-Roch theorem},
   note={Dedicated to William Fulton on the occasion of his 60th birthday},
   journal={Michigan Math. J.},
   volume={48},
   date={2000},
   pages={473--482},
   issn={0026-2285},
   doi={10.1307/mmj/1030132729},
}

\bib{Ooishi-CR}{article}{
   author={Ooishi, Akira},
   title={Castelnuovo's regularity of graded rings and modules},
   journal={Hiroshima Math. J.},
   volume={12},
   date={1982},
   number={3},
   pages={627--644},
   issn={0018-2079},
}

\bib{PR-ARR}{article}{
   author={Pink, Richard},
   author={R\"{o}ssler, Damian},
   title={On the Adams-Riemann-Roch theorem in positive characteristic},
   journal={Math. Z.},
   volume={270},
   date={2012},
   number={3-4},
   pages={1067--1076},
   issn={0025-5874},
   doi={10.1007/s00209-011-0841-7},
}

\bib{Quart-Loc}{article}{
   author={Quart, George},
   title={Localization theorem in $K$-theory for singular varieties},
   journal={Acta Math.},
   volume={143},
   date={1979},
   number={3-4},
   pages={213--217},
   issn={0001-5962},
   doi={10.1007/BF02392093},
}


\bib{Quillen-Higher}{article}{
   author={Quillen, Daniel},
   title={Higher algebraic $K$-theory: I [MR0338129]},
   conference={
      title={Cohomology of groups and algebraic $K$-theory},
   },
   book={
      series={Adv. Lect. Math. (ALM)},
      volume={12},
      publisher={Int. Press, Somerville, MA},
   },
   date={2010},
   pages={413--478},
}

\bib{Quillen-Homology}{article}{
   author={Quillen, Daniel},
   title={On the (co-) homology of commutative rings},
   conference={
      title={Applications of Categorical Algebra},
      address={Proc. Sympos. Pure Math., Vol. XVII, New York},
      date={1968},
   },
   book={
      publisher={Amer. Math. Soc., Providence, R.I.},
   },
   date={1970},
   pages={65--87},
}

\bib{Ro-Loc}{article}{
   author={R\"{o}ssler, Damian},
   title={A local refinement of the Adams-Riemann-Roch theorem in degree
   one},
   conference={
      title={Arithmetic L-functions and differential geometric methods},
   },
   book={
      series={Progr. Math.},
      volume={338},
      publisher={Birkh\"{a}user/Springer, Cham},
   },
   date={[2021] },
   pages={213--246},
}


\bib{Savin-MO}{misc}{
title={Answer to the question 'Local versus global embedding dimension' asked by R. Thomas.\\},
note={\url{https://mathoverflow.net/questions/351280/local-versus-global-embedding-dimension}}
}

\bib{StacksProject}{misc}{
title={The Stacks project},
note={\url{https://stacks.math.columbia.edu/}}
}

\bib{Thom-Alg}{article}{
   author={Thomason, R. W.},
   title={Algebraic $K$-theory of group scheme actions},
   conference={
      title={Algebraic topology and algebraic $K$-theory},
      address={Princeton, N.J.},
      date={1983},
   },
   book={
      series={Ann. of Math. Stud.},
      volume={113},
      publisher={Princeton Univ. Press, Princeton, NJ},
   },
   date={1987},
   pages={539--563},
}

\bib{Thom-L}{article}{
   author={Thomason, R. W.},
   title={Une formule de Lefschetz en $K$-th\'{e}orie \'{e}quivariante alg\'{e}brique},
   language={French},
   journal={Duke Math. J.},
   volume={68},
   date={1992},
   number={3},
   pages={447--462},
   issn={0012-7094},
   doi={10.1215/S0012-7094-92-06817-7},
}

\bib{TT}{article}{
   author={Thomason, R. W.},
   author={Trobaugh, Thomas},
   title={Higher algebraic $K$-theory of schemes and of derived categories},
   conference={
      title={The Grothendieck Festschrift, Vol. III},
   },
   book={
      series={Progr. Math.},
      volume={88},
      publisher={Birkh\"{a}user Boston, Boston, MA},
   },
   date={1990},
   pages={247--435},
}

\bib{Trung-C}{article}{
   author={Ng\^{o} Vi\^{e}t Trung},
   title={The Castelnuovo regularity of the Rees algebra and the associated
   graded ring},
   journal={Trans. Amer. Math. Soc.},
   volume={350},
   date={1998},
   number={7},
   pages={2813--2832},
   issn={0002-9947},
   doi={10.1090/S0002-9947-98-02198-9},
}


\bib{Washington-Cyc}{book}{
   author={Washington, Lawrence C.},
   title={Introduction to cyclotomic fields},
   series={Graduate Texts in Mathematics},
   volume={83},
   edition={2},
   publisher={Springer-Verlag, New York},
   date={1997},
   pages={xiv+487},
   isbn={0-387-94762-0},
   doi={10.1007/978-1-4612-1934-7},
}

\bib{WW-CMA}{book}{
   author={Whittaker, E. T.},
   author={Watson, G. N.},
   title={A course of modern analysis---an introduction to the general
   theory of infinite processes and of analytic functions with an account of
   the principal transcendental functions},
   edition={5},
   note={Edited by Victor H. Moll;
   With a foreword by S. J. Patterson;},
   publisher={Cambridge University Press, Cambridge},
   date={2021},
   pages={lii+668},
   isbn={978-1-316-51893-9},
}
\end{biblist}
\end{bibdiv}
\end{document}